\documentclass[a4paper, 10pt]{article}
\usepackage{subfiles}
\usepackage[toc,page]{appendix}
\usepackage[bottom]{footmisc}
\usepackage{microtype}
\usepackage{tocbibind}
\usepackage{stmaryrd}   
\usepackage{amsmath}
\usepackage{verbatim}
\usepackage{mathrsfs}
\usepackage{mathtools}
\usepackage{amsthm} 
\usepackage{amssymb} 
\usepackage{faktor}
\usepackage[colorlinks={true},linkcolor={blue},citecolor={blue}, urlcolor={blue} ]{hyperref}
\usepackage{braket}
\usepackage{transparent}
\usepackage{graphicx}
\usepackage{color}
\usepackage{rotating}
\usepackage{tikz}
\usepackage{tikz-cd}
\usepackage{authblk}
\usepackage[margin=2.5 cm]{geometry}
\usepackage{adjustbox}

\usetikzlibrary{intersections}
\usetikzlibrary{spy} 
\usetikzlibrary{decorations.pathmorphing}
\usetikzlibrary{patterns}
\usetikzlibrary{backgrounds}

\DeclareMathOperator{\divi}{div}

\DeclareMathOperator{\Div}{Div}

\DeclareMathOperator{\spec}{Spec}

\DeclareMathOperator{\Pic}{Pic}
\DeclareMathOperator{\len}{length}

\DeclareMathOperator{\mult}{mult}

\DeclareMathOperator{\id}{id}

\DeclareMathOperator{\Princ}{Princ}

\DeclareMathOperator{\Hom}{Hom}

\DeclareMathOperator{\supp}{supp}

\DeclareMathOperator{\Ob}{Ob}

\numberwithin{equation}{section}

\theoremstyle{plain}
\newtheorem{theorem}{Theorem}[section]
\newtheorem{lemma}[theorem]{Lemma}
\newtheorem{proposition}[theorem]{Proposition}

\theoremstyle{definition}
\newtheorem{definition}[theorem]{Definition}

\theoremstyle{remark}
\newtheorem{example}[theorem]{Example}
\newtheorem{remark}[theorem]{Remark}

\newcommand{\pctext}[2]{\text{\parbox{#1}{\centering #2}}}

\newcommand{\sbt}{\,\begin{picture}(-1,1)(-1,-3)\circle*{3}\end{picture}\ }

\newcommand{\catname}[1]{\mathbf{#1}}

\begin{document}
 \title{Explicit Deligne pairing}
\author{Paolo Dolce}
\date{}
\maketitle
\begin{abstract}
We give an explicit formula for the Deligne pairing for a proper and flat morphisms  $f:X\to S$ of schemes, in terms of the determinant of cohomology. The whole construction is justified by an analogy with the intersection theory on non-singular projective algebraic varieties.
\end{abstract}
\makeatletter
\@starttoc{toc}
\makeatother

\setcounter{section}{-1}

\section{Introduction}
The intersection pairing between  two  divisors on a projective non-singular surface is the unique bilinear and symmetric pairing with values in $\mathbb Z$ that satisfies some very natural properties: it counts the number of intersection points (when the divisor are normal crossing) and it is invariant if we ``move'' any of the two divisors in their linear class of equivalence. With the same philosophy, such a definition of intersection pairing between divisors can be extended naturally for projective non-singular varieties of any dimension: we list a number of natural properties and we find a unique multi-linear  symmetric pairing satisfying them. It turns out that this unique intersection pairing on algebraic varieties can be expressed explicitly in terms of the Euler-Poincare characteristics of (invertible sheaves associated to the) divisors. For example, for a surface over a field $k$ we have the well known formula:
\begin{equation}\label{eq_intro}
C.D=\chi_k(\mathscr O_X)-\chi_k(\mathscr O_X(C)^{-1})-\chi_k(\mathscr O_X(D)^{-1})+\chi_k(\mathscr O_X(C)^{-1}\otimes\mathscr O_X(D)^{-1}) 
\end{equation}
which is involved in the proof of the Riemann-Roch theorem for surfaces (see for example \cite{Beau}).

For a relative scheme $X\to S$, if we don't appeal to any ``compactification arguments'' of $X$ and $S$, there is in general no hope for finding a non-trivial reasonable intersection pairing for divisors which is invariant up linear equivalence. Let's see a simple example in the case of an arithmetic surface $X\to \spec \mathbb Z$: consider a prime $p\in\spec\mathbb Z$, then the fibre $X_{p}$ is a principal vertical divisor on $X$. Let $D$ be an effective, irreducible, horizontal divisor on $X$, then certainly $D$ meets $X_p$, so on one hand $D.X_p>0$ since our phantomatic intersection pairing should count the number of intersection points with multiplicity; but on the other hand we  said that $X_p$ is principal, which means $D.X_p=0$.

The closest object to an intersection pairing on a relative scheme $X\to S$ of relative dimension $n$ is the Deligne pairing. It is a map
$$\left<\;,\ldots,\,\right>_{X/S}:\operatorname{Vec}_1(X)^{n+1}\to \operatorname{Vec}_1(S)\,,$$ 
where $\operatorname{Vec}_1(\cdot)$ denotes the set of invertible sheaves, which descends to a symmetric, multi-linear map at the level of Picard groups. This pairing is of crucial importance in arithmetic geometry, since it gives ``the schematic contribution'' to the Arakelov intersection number.

The Deligne pairing was originally constructed  by Deligne in \cite{Del} for arithmetic surfaces and then generalised to any dimension in \cite{Elk}, \cite{Zh} and \cite{Gar}. Its definition was not built as the unique solution of a universal problem, but it was rather  constructed locally in terms of meromorphic sections of invertible sheaves. A set of axioms that uniquely identify the Deligne pairing have been found recently in the preprint \cite{Xia}.

For arithmetic surfaces one can show the following isomorphism of invertible sheaves which turns out to be crucial in the proof of Faltings-Riemann-Roch theorem (see for example \cite{Del} and \cite{MB} for more details):
\begin{equation}\label{eq_intro1}
\left<\mathscr L,\mathscr M\right>_{X/S}\cong \det Rf_\ast(\mathscr O_X)\otimes (\det Rf_\ast(\mathscr L))^{-1}\otimes(\det Rf_\ast(\mathscr M))^{-1}\otimes\det Rf_\ast(\mathscr L\otimes \mathscr M)\,.
\end{equation}
One can notice immediately the similarities between equations (\ref{eq_intro}) and (\ref{eq_intro1}). The only substantial difference is that for algebraic surfaces we use the Euler-Poincare characteristic, whereas for arithmetic surfaces we use the determinant of the cohomology. Such a distinction makes perfectly sense since the determinant of the cohomology is constructed to be the arithmetic analogue of the Euler-Poincare characteristic.

At this point the natural question is the following one: is it possible to give an explicit definition of the Deligne pairing (in the most general case) in terms of the determinant of cohomology\footnote{This was briefly conjectured already in \cite{Elk}: `` (...) Dans le cas g\'en\'eral, c'est-\`a-dire en dimension quelconque, et sans hypoth\`ese de lissit\' e, l'intersection doit aussi s'exprimer en termes de d\'eterminants d'images directes (...) malheureusement, pour l'instant, des probl\'e mes de signe obscurcissent s\'e rieusement la situation.''}? In this paper we give an affirmative answer. By working in complete analogy of the theory of algebraic varieties, we write down a simple explicit formula for the Deligne pairing in terms of the determinant of cohomology. Let $f:X\to S$ be a proper, flat morphism of integral Noetherian schemes, and assume that $f$ has pure dimension $n$, then we put
\begin{align}\label{eq_intro2}
\left\{\begin{array} {lll}
\left<\mathscr L_0,\ldots, \mathscr L_n\right>_{X/S}:={\det Rf_\ast\left(c_1(\mathscr L^{-1}_0)c_1(\mathscr L^{-1}_1)\ldots c_1(\mathscr L^{-1}_n)\mathscr O_X\right)}\, & \textrm{if } n>0\\
&\\
\left< \mathscr L_0\right>_{X/S}:={\det Rf_\ast\left(c_1(\mathscr L^{-1}_0)\mathscr O_X\right)}^{-1}=N_{X/S}(\mathscr L_0)& \textrm{if } n=0\\
\end{array}\right.
\end{align}
\noindent where $c_1(\mathscr L)\mathscr F:=\mathscr F-\mathscr L^{-1}\otimes \mathscr F\in K_0(X)$ for any coherent sheaf $\mathscr F$ and any invertible sheaf $\mathscr L$ (here only the class of $\mathscr F$ in the Grothendieck group matters). Moreover $N_{X/S}$ is the norm, relative to $f$, of an invertible sheaf.  We show that definition  (\ref{eq_intro2}) satisfies the axioms of  \cite{Xia}, and this implies that our definition is exactly the Deligne pairing.

Let's mention some other papers that previously investigated in our direction: an explicit formula for the Deligne pairing  when $X$ and $S$ are integral schemes over $\mathbb C$ was announced in \cite{BSW}, although a complete proof is not given. The approach of \cite{BSW} is essentially different from ours, indeed the authors work on local trivializations of invertible sheaves. A  more complicated expression of the Deligne pairing in terms of symmetric difference of the  functor $\det Rf_\ast$ is proved in \cite{Ducr} with heavy usage of category theory (see also \cite[Appendix A]{bouck}). Moreover, when $\mathscr L$ is very ample on the fibres, an explicit expression of the Deligne pairing $\left<\mathscr L,\ldots,\mathscr L\right>$,  is given in \cite{Phon} as the leading term of the Knudsen-Mumford expansion of $\det Rf_\ast\left(\mathscr L^k\right)$. \\

This paper is organized in the following way: in section \ref{endo} we introduce the map $c_1(\mathscr L):K_0(X)\to K_0(X)$  with all its properties. Section \ref{int_var} is a review of intersection theory for algebraic varieties and it gives to the reader the philosophical guidelines for the case of relative schemes. In section \ref{rel_sch} we give the axioms of the Deligne pairing and we show with all details that if such pairing exists, then it must be unique (we follow \cite{Xia}). Afterwards we show that  the pairing (\ref{eq_intro2}) satisfies all the axioms. appendix \ref{det} is a review of the determinant of cohomology, this part is crucial in order to understand section \ref{rel_sch}.  Finally, in appendix \ref{ori_constr} we review in all details the original construction of Deligne pairing of \cite{Del} (very often this construction is just sketched in the literature).

\paragraph{Acknowledgements.}
The author wants to express his gratitude to R.S. De Jong for his time spent in discussing the topic during summer 2019 in Nottingham and for his precious insights. A special thanks goes also to P. Corvaja, I. Fesenko, S. Urbinati, F. Zucconi and to the anonymous referee for his/her careful reading and for providing toughtful comments.\\

\noindent This research was supported by the Italian national grant ``Ing. Giorgio Schirillo" conferred by INdAM and partially by the EPSRC programme grant EP/M024830/1 (Symmetries and correspondences: intra-disciplinary developments and applications).

\section{An endomomorphism of the group $K_0(X)$}\label{endo}
 Let's briefly recall the abstract construction of the Groethendieck group $K_0(\catname C)$. Fix an abelian category  $\catname{C}$ and let $F(\catname C)$ be the free abelian group over the set $\Ob(\catname C)/\cong$, where $\cong$ is the isomorphism relation. If $C\in \Ob(\catname C)$, then $(C)$ denotes isomorphism class in $\Ob(\catname C)/\cong$. To any  short exact sequence in $\catname C$:
$$\mathcal S\colon \quad 0\to C'\to C\to C''\to 0$$
we associate an element $ Q(\mathcal S):=(C)- (C')-(C'')\in F(\catname C)$. Now $H(\catname C)$ is the subgroup of $F(\catname C)$ generated by all the elements $Q(\mathcal S)$ for $\mathcal S$ running over all short exact sequences. Then:
$$K_0(\catname C):= F(\catname C)/H(\catname C)\,,$$
and $[C]\in K_0(\catname C)$ denotes the equivalence class associated to $C\in \Ob(\catname C)$.

Let's fix a Noetherian scheme $X$, then $K_0(X):=K_0(\catname {Coh}(X))$, where $\catname {Coh}(X)$ is the category of coherent sheaves on $X$. From now on, by an abuse of notation we identify any coherent sheaf $\mathscr F$ with its class in $K_0(X)$. In this paper, with the notation $\catname {Coh}_r(X)$ we denote the category of coherent sheaves on $X$ whose support has dimension at most $r$, and we define $K_{0,r}(X):= K_0(\catname {Coh}_r(X))$. Clearly when $0\le i\le j$, then $K_{0,i}(X)\subseteq K_{0,j}(X)$.

For any invertible sheaf $\mathscr L$ on $X$ we define a map:

\begin{equation}
\begin{aligned}
 c_1(\mathscr L):K_0(X)&\to  K_0(X)\\
\mathscr F &\mapsto  c_1(\mathscr L)\mathscr F:= \mathscr F-\mathscr L^{-1}\otimes \mathscr F\,.
\end{aligned}
\end{equation}
Note that it is well defined because tensoring with an invertible sheaf is an exact functor, moreover it defines and endomorphism of the group $K_0(X)$. Since the notation for the function $c_1(\mathscr L)$ is multiplicative, the symbol $c_1(\mathscr L)c_1(\mathscr L')$ denotes the composition of functions. The properties of the operator $c_1(\mathscr L)$ are well described in \cite[Appendix B]{Kl}, so here we just recall them.

\begin{proposition}\label{propc_1}
The following properties hold for the operator $c_1(\mathscr L)$:
\begin{itemize}
\item[$(i)$] $c_1(\mathscr L)c_1(\mathscr M)=c_1(\mathscr L)+ c_1(\mathscr M)-c_1(\mathscr L\otimes \mathscr M)$, where clearly the sum is taken in $\operatorname{End}(K_0(X))$.
\item[$(ii)$] $c_1(\mathscr M) c_1(\mathscr L)=c_1(\mathscr L) c_1(\mathscr M).$
\item[$(iii)$] If $Z\subset X$ is a closed subscheme  and $\mathscr L_{|Z}=\mathscr O_Z(D)$ where $D$ is an effective Cartier divisor on $Z$, then $c_1(\mathscr L) \mathscr O_Z=\mathscr O_D$.
\end{itemize}
\end{proposition}
\proof
Both sides of the equality in $(i)$ applied to $\mathscr F$ expand to:

\begin{equation}\label{e_comm}
\mathscr F-\mathscr L^{-1}\otimes \mathscr F-\mathscr M^{-1}\otimes \mathscr F\mathscr +\mathscr L^{-1}\otimes\mathscr M^{-1}\otimes\mathscr F\,.
\end{equation}
$(ii)$ Follows easily by looking at equation (\ref{e_comm}). For $(iii)$ Consider the short exact sequence
$$0\to\mathscr O_Z(-D)\to\mathscr O_Z\to\mathscr O_D\to 0\,.$$
\endproof
\begin{proposition}
Let $\mathscr F\in K_{0,r}(X)$ and let $Z_1,\ldots, Z_s$ the $r$-dimensional irreducible components of $\supp(\mathscr F)$ whose  generic points are denoted respectively by $z_i$. Let $n_i=\len \mathscr F_{z_i}$. Then in $K_{0,r}(X)$ we have the equality:
$$\mathscr  F\equiv \sum^s_{i=1} n_i\mathscr O_{Z_i}\mod K_{0,r-1}(X)\,.$$
\end{proposition}
\proof
See \cite[Lemma B4]{Kl}.
\endproof

\begin{proposition}\label{decreasing_sup}
let $\mathscr L$ an invertible sheaf on $X$, then $c_1(\mathscr L)K_{0,r}(X)\subset K_{0,r-1}(X)$ for any $r\ge 0$.
\end{proposition}
\proof
See \cite[Lemma B5]{Kl}.
\endproof

\begin{remark}
The operator $c_1(\mathscr L)$ can be ``extended'' to bounded complexes of coherent sheaves on $X$. Let $\mathscr F^\bullet$ be a bounded complex of objects in $\catname{Coh}(X)$ then we can define:
$$c_1(\mathscr L)\mathscr F^\bullet:=\sum_i(-1)^ic_1(\mathscr L)\mathscr F^i\in K_0(X)\,.$$
Such a map is clearly zero on short exact sequences. 
\end{remark}

\section{Intersection theory for algebraic varieties}\label{int_var}

\begin{definition}\label{inters_var}
Let $X$ be a $n$-dimensional projective, non-singular algebraic variety over a field $k$. An \emph{intersection pairing} on $X$ is a map:
\begin{equation}
\begin{aligned}
\Div(X)^n&\to  \mathbb Z\\
(D_1,\ldots,D_n) &\mapsto  D_1.D_2.\,\ldots\,.D_n
\end{aligned}
\end{equation}
satisfying the following properties:
\begin{itemize}
\item[$(1)$] It is symmetric and $\mathbb Z$-multilinear.
\item[$(2)$] It descends to a pairing $\Pic(X)^n\to  \mathbb Z$.
\item[$(3)$] Let $D_i$ be a prime divisor for any $i$  and let $e_{i,x}\in\mathscr O_{X,x}$ be a local equation of $D_i$ at the point $x$. Assume that for all $x$ in the support of all divisors $D_i$, the $e_{i,x}$'s form a regular sequence in  $\mathscr O_{X,x}$ (i.e. the divisors are in general position), then: 
$$D_1.D_2.\,\ldots\,.D_n=\sum_{x\in\cap D_i}\len\frac{\mathscr O_{X,x}}{(e_{1,x}, e_{2,x},\ldots, e_{n,x})}$$  
\end{itemize}
\end{definition}
Now we show that if an intersection pairing exists, it is uniquely defined by the three axioms of Definition \ref{inters_var}.
\begin{proposition}
If an intersection pairing exists, then it is unique. 
\end{proposition}
\proof
Let $\left<\;\right>_1$ and  $\left<\;\right>_2$ two pairings satisfying the axioms $(1)-(3)$ and fix $D_1,\ldots,D_n\in \Div(X)^n$;  by $(1)$ we can assume that all $D_i$ are prime. Thanks to Chow's moving lemma we can find some divisors $D'_i$ such that $D\sim D_i'$ and $D'_1,\ldots,D'_n$ are in general position. Therefore by using $(2)$ and $(3)$ we get:
$$\left<D_1,\ldots, D_n\right>_1=\left<D'_1,\ldots, D'_n\right>_1=\left<D'_1,\ldots, D'_n\right>_2=\left<D_1,\ldots, D_n\right>_2\,.$$
\endproof
The remaining part of this section is devoted to give the explicit expression of the intersection pairing on $X$ as in \cite{Sn1}, \cite{Sn2} and later \cite{Ca}; then we see that the axioms of definition  \ref{inters_var} are satisfied. Such an intersection pairing uses the endomorphism defined in section \ref{endo} and the Euler-Poincare characteristic for coherent sheaves. 

We actually give a definition of the intersection pairing in a more general setting, in fact we will assume that $X$ is a relative scheme over a scheme $S$, and we define a ``partial''  intersection number for a particular subclass of divisors.

From now on, in this section we assume that $X\to S$ is a flat and proper morphism of integral Noetherian schemes. Let's denote with $\catname{Coh}(X/S)$ the category of coherent sheaves on $X$ whose schematic support is proper  over a $0$-dimensional subscheme of $S$. Moreover $\catname{Coh}_r(X/S)$ is the subcategory of  $\catname{Coh}(X/S)$ made of sheaves whose support has dimension at most $r$. The motivation behind the restriction to sheaves with this kind of support is that for any $\mathscr F\in \catname{Coh}(X/S)$ we have a well defined notion of Euler-Poincare characteristic. In fact if $T$ is the schematic support of $\mathscr F$ and $S_0=f(T)$, we know that $S_0$ is Noetherian of dimension $0$, so $S_0=\spec A$ with $A$  artinian; at this point we can put

$$\chi_S(\mathscr F):=\sum_{i\ge 0} (-1)^i\len_A H^i(X, \mathscr F)\,.$$
When $S=\spec k$, then $\chi_S$ is the usual Euler-Poicare characteristic (for coherent sheaves with proper support). Thanks to the ``additivity'' of $\chi_S$ with respect to short exact sequences, it is immediate to notice that we have a naturally induced group homomorphism $\chi_S:K_0(\catname{Coh}_r(X/S))\to\mathbb Z$.
\begin{definition}
Let $X\to S$ be as above and consider $\mathscr F\in\catname{Coh}_r(X/S)$.  Then the \emph{intersection number of the invertible sheaves $\mathscr L_1,\ldots\mathscr L_r$ (with respect to $\mathscr F$)} is defined as:

$$(\mathscr L_1.\mathscr L_2.\,\ldots\,.\mathscr L_r,\mathscr F):=\chi_S\left(c_1(\mathscr L_1)c_1(\mathscr L_2)\ldots c_1(\mathscr L_r)\mathscr F\right)\,.$$
When $\mathscr F=\mathscr O_X$, which implies $r\ge\dim(X)$, we put for simplicity 

$$\mathscr L_1.\mathscr L_2.\,\ldots\,.\mathscr L_r:=(\mathscr L_1.\mathscr L_2.\,\ldots\,.\mathscr L_r,\mathscr O_X)\,.$$
Moreover if $\mathscr L_i=\mathscr O_X(D_i)$ for a Cartier divisor $D_i$ on $X$, then:
$$D_1.D_2.\,\ldots\,.D_r:=\mathscr  O_X(D_1).\mathscr  O_X(D_2).\,\ldots\,.\mathscr  O_X(D_r)\,.$$   

\end{definition}

\begin{example}
If $X$ is a surface over $k$ and $C,D$ are two divisors, then:
$$C.D=\chi_k(c_1(\mathscr O_X(C)) c_1(\mathscr O_X(D))\mathscr O_X)=\chi_k(c_1(\mathscr O_X(C))(\mathscr O_X-\mathscr O_X(D)^{-1}))=$$
$$=\chi_k(\mathscr O_X)-\chi_k(\mathscr O_X(C)^{-1})-\chi_k(\mathscr O_X(D)^{-1})+\chi_k(\mathscr O_X(C)^{-1}\otimes\mathscr O_X(D)^{-1}) \,.$$
\end{example}

The mere definitions tell us that we can intersect a number of divisors which is greater or equal to the dimension on $X$. On the other hand the next lemma shows that the intersection of a number of divisor which is strictly bigger than the dimension of $X$ is always $0$.
\begin{lemma}\label{no_more_dim}
If $\mathscr F\in\catname{Coh}_{r}(X/S)$, then $(\mathscr L_1.\mathscr L_2.\,\ldots\,.\mathscr L_{r+1},\mathscr F)=0$.
\end{lemma}
\proof
It follows directly from Proposition \ref{decreasing_sup}.
\endproof
\begin{proposition}\label{bilinear}
The intersection number of $\mathscr L_1,\ldots,\mathscr L_m$ with respect to $\mathscr F$ is a $\mathbb Z$-multilinear map in the $\mathscr L_j$'s (the operation is the tensor product).
\end{proposition}
\proof
Follows by Proposition \ref{propc_1}(i) and Lemma \ref{no_more_dim}.
\endproof
\begin{proposition}\label{mor_int}
Let $g:X'\to X$ be an morphism of $S$-schemes and let $\mathscr F\in\catname{Coh}_r(X'/S)$, then
$$(g^\ast\mathscr L_1.g^\ast\mathscr L_2\,.\ldots g^\ast\mathscr L_n, \mathscr F)=(\mathscr L_1.\mathscr L_2\,.\ldots \mathscr L_n, g_\ast\mathscr F)\,.$$ 
\end{proposition}
\proof
See \cite[Lemma B.15]{Kl}.
\endproof

We can give the explicit expression of the intersection number on varieties:
\begin{proposition}
Let $X$ be a non-singular algebraic variety of dimension $n$ over a field $k$. The pairing $$(D_1,\ldots,D_n)\mapsto D_1.D_2.\,\ldots\,.D_n$$
defines the intersection number on $X$.
\end{proposition}
\proof
Axiom $(1)$ is satisfied thanks to Proposition \ref{bilinear}. Axiom $(2)$ is obvious and axiom $(3)$ is \cite[IV, Theorem 2.8]{Kol}.
\endproof

Finally we state a proposition regarding the intersection along fibres:

\begin{proposition}\label{gen_fiber}
Let $s\in S$ be a closed point and let $X_s$ be the fibre over $b$. Assume that $\dim(X_s)=d$, then the map:
$$s\mapsto (\mathscr L_1,\ldots,\mathscr L_d; \mathscr O_{X_s})$$
is locally constant on $S$.
\end{proposition}
\proof
See \cite[VI, Proposition 2.10]{Kol}.
\endproof

\section{The case of schemes over a general base}\label{rel_sch}

\subsection{Multi-monoidal and symmetric functors}
The Deligne pairing will be expressed as a collection of functors, so in this section we recall what is the functorial equivalent of a multi-linear homomorphism of abelian groups. 

We assume that the reader is familiar with some basic notions of category theory and the concept of \emph{Picard groupoid}. Roughly speaking a  Picard groupoid is a category where the morphisms are all invertible and moreover there is a ``group-like'' operation between the object of the category. A simple example is the \emph{Picard category} $\catname{Pic}(X)$, made of all invertible sheaves on a scheme $X$, and where the morphisms are just the isomorphisms. The ``operation'' in $\catname{Pic}(X)$  is clearly the tensor product of invertible sheaves and the identity element is the structure sheaf. The morphisms we want to consider between Picard groupoids are monoidal functors, i.e. functors that preserve the monoidal structure of the categories:

For the remaining part of this subsection we fix two Picard groupoids $(\catname{C},\otimes)$ and $(\catname{D},\otimes)$.
\begin{definition}
A monoidal functor $\catname{C}\to\catname{D}$ is a collection $(F,\epsilon, \mu)$ where $\mu:=\{\mu_{X,Y}\}_{X,Y\in\text{Obj(}\catname C)}$, satisfying the following properties:
\begin{enumerate}
\item[\sbt] $F:\catname{C}\to\catname{D}$ is a functor.
\item[\sbt] $\epsilon\colon F(1_{\catname C})\xrightarrow{\cong} 1_{\catname D}$ is an isomorphism.
\item[\sbt] $\mu_{X,Y}\colon F(X)\otimes F(Y)\xrightarrow{\cong} F(X\otimes Y)$ is an isomorphism functorial in $X$ and $Y$ which satisfies associativity and unitality in the obvious categorical sense.
\end{enumerate}
For simplicity we often omit $\epsilon$ and $\mu$ and we say that $F$ is a monoidal functor between  $\catname C$ and $\catname{D}$. In symbols we write $F\in L^1(\catname C,\catname D)$.
\end{definition}

\begin{definition}\label{nat1}
A \emph{natural transformation} between monoidal functors $(F, \epsilon,\mu)$ and $(F',\epsilon', \mu')$ is a monoidal natural transformation $\alpha:F\to F'$ which maps  $\epsilon$ to $\epsilon'$ and $\mu$ to $\mu'$.
\end{definition}
In order to give the next definition we need to introduce some notations. An object of the category $\catname{C}^n$ (i.e. a $n$-uple of objects of $\catname{C}$) is denoted by $ X=(X_1,\ldots, X_n)$. Let $X, Y\in \catname{C}^n$  and let $i\in\{1,\ldots,n\}$ such that for any $j\in\{1,\ldots,n\}$ with $j\neq i$ we have $X_j=Y_j$, then we define $ X\otimes_i Y\in\catname{C}^n$ in the following way:

$$
(X\otimes_i Y)_j=\left\{\begin{array} {ll}
X_j\otimes Y_j & \text{if } i=j\\
X_j & \text{if } i\neq j\\
\end{array}\right.
$$

\begin{definition}
A multi-monoidal\footnote{Very often in literature one can find the term multi-additive.} functor  $\catname{C}^n\to\catname{D}$ is the datum of
\begin{enumerate}
\item[\sbt] A functor $F:\catname{C}^n\to \catname{D}$.
\item[\sbt] For any functor  $F':\catname{C}\to\catname{D}$ obtained by fixing $n-1$ components in $\catname C$, we have a collection $\mu'$ such that $(F',\mu')$ is a monoidal functor $\catname{C}$ and $\catname D$.
\item[\sbt] For every $i,j\in\{1,\ldots, n\} $  and $X,Y,Z,W\in\catname{C}^n$ such that   $X_k=Y_k=Z_k=W_k$ for all $k\neq i,j$, a  commutative diagram:
$$
\adjustbox{scale=0.8,center}{%
\begin{tikzcd}
 & F(X)\otimes F(Y)\otimes  F(Z)\otimes F(W) \arrow{dl}\arrow{dr} &\\
(F(X)\otimes_j F(Y))\otimes (F(Z)\otimes_j F(W))\arrow{d}& & (F(X)\otimes_i F(Y))\otimes (F(Z)\otimes_i F(W))\arrow{d}\\
F((X\otimes_j Y)\otimes_i (Z\otimes_j W))\arrow["="]{rr}&& (F(X)\otimes_i F(Y))\otimes_j (F(Z)\otimes_i F(W))
\end{tikzcd}
}
$$
\end{enumerate}
\end{definition}
The notion of symmetry is what one expects:
\begin{definition}
A multi-monoidal functor  $\catname{C}^n\to\catname{D}$ is \emph{symmetric} if for any $c_i\in\catname C$ and any permutation $\sigma\in\Sigma_n$ we have $F(c_1,\ldots c_n)\cong F(c_{\sigma(1)},\ldots c_{\sigma(n)})$.
\end{definition}

The set of symmetric multi-monoidal functors from $\catname{C}^n$ to $\catname D$ is denoted by $L^n(\catname C,\catname D)$.

\begin{definition}
\emph{A natural transformation between two multi-monoidal functors} $F,F'\in L^n(\catname C,\catname D)$ is a functorial isomorphism $\alpha: F\to F'$ which restricts to a natural transformation to each component in the sense of definition \ref{nat1}.
\end{definition}

\subsection{Axiomatic Deligne pairing}
The Deligne pairing was introduced in \cite{Del} as a bilinear and symmetric map $\left<\,,\,\right>:\operatorname{Vec}_1(X)\times\operatorname{Vec}_1(X)\to\operatorname{Vec}_1(S)$, where $X\to S$ is an arithmetic surface. Such a definition requires the choice of meromorphic sections ``behaving well" on an open set, and  then clearly one has to show the independence with respect to this choice. The Deligne pairing satisfies some compatibility conditions with respect to the base change, the pullback functor and the norm functor. In \cite{Elk} Deligne's construction was extended straight away for proper  flat morphisms of integral schemes of any dimension.

Let $f:X\to S$ a proper flat morphism between Noetherian integral schemes, the guiding idea of this paper is that the  Deligne pairing relative to $f$ should be the generalisation the intersection pairing described in section \ref{int_var}.  We want to work in complete analogy with the case of algebraic varieties, so in this section we give a set of ``natural axioms'' that uniquely define the Deligne pairing\footnote{We follow \cite{Xia}, but we prefer to give a self contained presentation with all details.}. The explicit construction of the Deligne pairing will be carried out  in section  \ref{det_del}.

Let $X$ and $S$ be two Noetherian integral schemes, with the symbol $\mathcal F^n(X,S)$ we denote the set of all proper flat morphisms $X\to S$ of pure dimension $n$.

\begin{definition}\label{axiom_deligne}
A \emph{Deligne pairing} consists of the following data for any $f\in\mathcal F^n(X,S)$ where $X$ and $S$ are two Noetherian integral schemes: a functor 
$$\left< \,\;,\;,\ldots,\;\right>_f=\left< \,\;,\;,\ldots,\;\right>_{X/S}\in L^{n+1}(\catname{Pic}(X), \catname{Pic}(S))$$
and a collection of natural transformations $\alpha,\beta,\gamma,\delta$ described below:
\begin{enumerate}
\item[(1)] For any commutative square given by a base change $g:S'\to S$ which is proper, flat and with connected fibres 

$$
\begin{tikzcd}[row sep=large, column sep = huge]
X'=X\times_S S'\arrow["f'"]{r}\arrow["g'"]{d} & S'\arrow["g"]{d}\\
X\arrow["f"]{r} & S
\end{tikzcd}
$$
a natural transformation  $\alpha_{f,g}$ between multi-monoidal functors $\catname{Pic}(X)^{n+1}\to \catname{Pic}(S')$ such that
 $$\alpha_{f,g}:g^\ast\left<\mathscr L_0,\ldots,\mathscr L_n\right>_{X/S}\xrightarrow{\cong}\left<g'^\ast \mathscr L_0,\ldots,g'^\ast\mathscr L_n\right>_{X'/S'}\,.$$
 
 \item[(2)] When $n>0$ and $D\in\Div(X)$ is an effective relative Cartier divisor, a natural transformation $\beta_{f,D}$ between multi-monoidal functors $\catname{Pic}(X)^{n}\to \catname{Pic}(S)$ such that
  $$\beta_{f,D}:\left<\mathscr L_1,\ldots,\mathscr L_n, \mathscr O_X(D)\right>_{X/S}\xrightarrow{\cong}\left<\mathscr L_1|_D,\ldots,\mathscr L_n|_D\right>_{D/S}\,.$$
Moreover $\beta_{f,D}$ is natural with respect to  base change in the following sense: for a  base change diagram as in axiom $(1)$ we have a commutative diagram:  
 $$
\begin{tikzcd}[row sep=large, column sep = huge]
\left<g'^\ast\mathscr L_1,\ldots,g'^\ast\mathscr L_n, g'^\ast\mathscr O_X(D)\right>_{X'/S'}\arrow["\beta_{f',g'^\ast D}"]{r}\arrow["\cong"]{d} & \left<g'^\ast\mathscr L_1|_{g'^\ast D},\ldots,g'^\ast\mathscr L_n|_{g'^\ast D}\right>_{{g'^\ast D}/S'}\arrow["\cong"]{d}\\
g^\ast\left<\mathscr L_1,\ldots,\mathscr L_n, \mathscr O_X(D)\right>_{X/S} & g^\ast \left<\mathscr L_1|_D,\ldots,\mathscr L_n|_D\right>_{D/S}\arrow["g^\ast \beta_{f,D}"]{l} 
\end{tikzcd}
$$
Where the vertical isomorphisms are given by $\alpha_{f,g}$ (remember that $g'^\ast\mathscr O_X(D)=\mathscr O_{X'}(g'^\ast D)$).
 \item[(3)] When $n>0$, a natural transformation $\gamma_{f}$ between multi-monoidal functors $\catname{Pic}(S)\times\catname{Pic}(X)^{n}\to \catname{Pic}(S)$
 such that
  $$\gamma_{f}:\left<f^\ast\mathscr L,\mathscr L_1\ldots,\mathscr L_n\right>_{X/S}\xrightarrow{\cong}\mathscr L^{(\mathscr L_1|_{X_s}.\mathscr L_2|_{X_s}.\ldots\mathscr L_n|_{X_s};\mathscr O_{X_s})}$$
where $X_s$ is a generic fibre of $f$ (see Proposition \ref{gen_fiber}). Moreover $\gamma_{f}$ is natural with respect to base change in the following sense: for a  base change diagram as in axiom $(1)$ we have a commutative diagram:  
 $$
\begin{tikzcd}[row sep=large, column sep = huge]
\left<g'^\ast f^\ast\mathscr L,g'^\ast\mathscr L_1\ldots,g'^\ast\mathscr L_n\right>_{X'/S'}\arrow["\gamma_{f'}"]{r}\arrow["\cong"]{d} &  g^\ast \mathscr L^{(g'^\ast \mathscr L_1|_{g'^\ast X_s}. g'^\ast \mathscr L_2|_{g'^\ast X_s}.\ldots g'^\ast\mathscr L_n|_{g'^\ast X_s};\mathscr O_{g'^\ast X_s})}\arrow["="]{d}\\
g^\ast\left<f^\ast\mathscr L,\mathscr L_1\ldots,\mathscr L_n\right>_{X/S} & g^\ast\mathscr L^{(\mathscr L_1|_{X_s}.\mathscr L_2|_{X_s}.\ldots\mathscr L_n|_{X_s};\mathscr O_{X_s})}\arrow["g^\ast \gamma_{f}"]{l} 
\end{tikzcd}
$$
where the vertical isomorphism is given $\alpha_{f,g}$  and the equality follows from proposition \ref{mor_int} and the properties of $g$.

\item[(4)] When $n=0$, a natural transformation $\delta_{f}$ between monoidal functors $\catname{Pic}(X)\to \catname{Pic}(S)$
such that
  $$\delta_{f}:\left<\mathscr L\right>_{X/S}\xrightarrow{\cong} N_{X/S}(\mathscr L)$$
 where $N_{X/S}$ is the norm of $f$ (see Definition \ref{norm}). Moreover $\delta_f$ is natural with respect to base change in the following sense: for  a base change diagram as in axiom $(1)$ we have a commutative diagram:
$$
\begin{tikzcd}[row sep=large, column sep = huge]
\left<g'^\ast\mathscr L\right>_{X'/S'}\arrow["\delta_{f'}"]{r}\arrow["\cong"]{d} & N_{X'/S'}(g'^\ast\mathscr L)\arrow["\cong"]{d}\\
g^\ast\left<\mathscr L\right>_{X/S} & g^\ast N_{X/S}(\mathscr L)\arrow["g^\ast \delta_f"]{l} 
\end{tikzcd}
$$
where the vertical isomorphisms are given respectively by $\alpha_{f,g}$ and  thanks to the properties of the norm.
\end{enumerate} 
\end{definition}
We have to show that if a Deligne pairing exists, then it is unique. Roughly speaking  we will show that any two pairings $(\left< \,\;,\;,\ldots,\;\right>^i,\alpha^i,\beta^i,\gamma^i,\delta^i)$, with $i=1,2$, satisfying the axioms of definition \ref{axiom_deligne} are related by  natural transformation of functors that respects all the data. We will work by induction on the relative dimension of the morphism $f$. Note that we cannot use straight away property $(2)$ to  pass from relative dimension $n$ to $n-1$, since the whole construction would depend on the choice of a relative divisor $D$, whereas we want our constructions to be natural in a functorial way. So, let's describe a general well known procedure to reduce the relative dimension of $f$ by using a canonical choice of a relative Cartier divisor. It is called \emph{universal extension}.

Let $f\in\mathcal F^n(X,S)$ and let $\mathscr L$ be an invertible sheaf on $X$. We assume that $\mathscr L$ is \emph{sufficiently ample with respect to $f$}, i.e. that the following properties are satisfied: $\mathscr L$ is very ample with respect to $f$ and $R^i f_\ast\mathscr L=0$ for $i>0$. 
\begin{remark}
The following properties hold for sufficient ampleness:
\begin{itemize}
\item[\sbt] It is preserved after base change.
\item[\sbt] If $f\in\mathcal F^n(X,S)$ and $\mathscr L$ is sufficiently ample on $X$, then $f_\ast \mathscr L$ is a locally free sheaf on $S$.
\item[\sbt] If $\mathscr L_0$ is an invertible sheaf on $X$, then there exists a sufficiently ample $\mathscr L$ such that $\mathscr L_0\otimes \mathscr L$ is sufficiently ample. In particular we can always find on $X$ a sufficiently ample invertible sheaf.
\end{itemize}
\end{remark}
Put $\mathscr M=(f_\ast\mathscr L)^{\vee}$ and let $\mathbb P:=\mathbb P_S(\mathscr M)$ be the projective vector bundle associated to $\mathscr M$, over $S$. Then we obtain the following base change diagram: 
\begin{equation}\label{diag_univ}
\begin{tikzcd}[row sep=large, column sep = huge]
\mathbb X:=\mathbb P\times_S X\arrow["p_2"]{r}\arrow["p_1"]{d} & X\arrow["f"]{d}\\
\mathbb P\arrow["\pi"]{r} & S
\end{tikzcd}
\end{equation}
Consider now the invertible sheaf $\mathscr L_f:=p_1^\ast\mathscr O_\mathbb P(1)\otimes p_2^\ast\mathscr L$ on $\mathbb X$. We want to construct a \emph{canonical global section} $\sigma$ of $\mathscr L_f$. It is enough to find a canonical non-zero element in $\mathscr L^{-1}_f$, because if $\phi\in\Hom(\mathscr O_X,\mathscr L_f)=\mathscr L^{-1}_f$ then we put $\sigma:=\phi_X(1)$.
First of all we construct a surjective canonical morphism 
$$\Psi:f^\ast \mathscr M\to\mathscr L \,.$$ 
Thanks to the properties of the pullback we have a canonical isomorphism $f^\ast \mathscr M \cong (f^\ast f_\ast \mathscr L)^\vee$. Since $\mathscr L$ is sufficiently ample, we have a canonical isomorphism   $(f^\ast f_\ast \mathscr L)^\vee\cong \mathscr L^\vee$. Moreover there is a surjective canonical map $\mathscr L^\vee\to\mathscr L$ given in the following way:
\begin{eqnarray*}
\Hom(\mathscr O_X(U), \mathscr L(U))& \to& \mathscr L(U)\\
\varphi_U &\mapsto& \varphi_U(1)\,.
\end{eqnarray*}
By taking all compositions, we finally get our surjective $\Psi$. We have to prove that $\Psi$ induces a canonical element in $\mathscr L^{-1}_f$ (in order to get $\sigma$). Note that $\mathscr L^{-1}_f$ is canonically isomorphic to $\Hom(\mathscr L^{-1}, {p_2}_\ast p_1^\ast \mathscr O_{\mathbb P}(1))$, but  
$${p_2}_\ast p_1^\ast \mathscr O_{\mathbb P}(1)=f^\ast \pi_\ast\mathscr O_{\mathbb P}(1)=f^\ast (\mathscr M^\vee)=(f^\ast \mathscr M)^\vee\,.$$
We conclude that the dual map of $\Psi$ induces the non-zero element of $\mathscr L^{-1}_f$ that we were searching for.

From now on we will say that the section $\sigma$ constructed above is the universal section relative to $\mathscr L$.  The following remark explains why we can use the universal section for our inductive step in the proof of uniqueness:

\begin{remark}
In \cite[2.2]{Gar} it is shown that $\sigma$ is a regular section, which is equivalent to say that the zero locus of $\sigma$ (considered with its reduced scheme structure):
$$Z(\sigma):=\{x\in\mathbb X\colon 0=\sigma(x)\in \mathscr L_f/\mathfrak m_x \mathscr L_f\}$$
is a relative Cartier divisor on $\mathbb X$. In this case we also have that $\mathscr L_f$ is canonically isomorphic to $\mathscr O_{\mathbb X}(Z(\sigma))$. Now consider the restriction 
$$p:=(p_1)|_{Z(\sigma)}: Z(\sigma )\to \mathbb P\,;$$
Let $U$ be the flat locus of $p$ and put $V:=p(U)$. Then $V$ is open in  $\mathbb P$, and we denote its closed complement with $W$, then we conclude that  
\begin{equation}\label{section_map}
p: Z(\sigma)-p^{-1}(W)\to V
\end{equation}
is flat of relative dimension $n-1$.
\end{remark}

The following theorem ensures the unicity of the Deligne pairing:
\begin{theorem}
The Deligne pairing is unique: given two sets of data $(\left< \,\;,\;,\ldots,\;\right>^i,\alpha^i,\beta^i,\gamma^i,\delta^i)$, with $i=1,2$, satisfying the conditions of Definition \ref{axiom_deligne}, there is a unique multi-monoidal morphism  $\left< \,\;,\;,\ldots,\;\right>^1\to \left< \,\;,\;,\ldots,\;\right>^2$ that transforms $\alpha^i,\beta^i,\gamma^i,\delta^i$ accordingly.
\end{theorem}
\proof
We proceed by induction on $n$. When $n=0$, the claim follows directly from axiom $(4)$. Let's work now with  $n>0$; first of all we want a functorial isomorphism:
\begin{equation}\label{iso_to_show}
\Psi(\mathscr L_0,\ldots, \mathscr L_n):\left<\mathscr L_0,\ldots,\mathscr L_n\right>^1_{X/S}\xrightarrow{\cong}\left<\mathscr L_0,\ldots,\mathscr L_n\right>^2_{X/S}\,.
\end{equation}
Let's first construct it by assuming that one invertible sheaf  $\mathscr L=\mathscr L_0$, is chosen sufficiently ample; we will denote it as $\Psi'(\mathscr L_0,\ldots, \mathscr L_n)$. Let's construct for $\mathscr L$ the base change diagram (\ref{diag_univ}), with the same notations. Then $\sigma$ is the universal section of $\mathscr L_f$ and we also have the map $p$ described in equation (\ref{section_map}). Thanks to \cite[Lemme 21.13.2]{EGAIV}, in order to give isomorphism (\ref{iso_to_show}), it is enough to give a functorial isomorphism:
\begin{equation}\label{iso_to_show1}
(\pi|_{V})^\ast\left<\mathscr L, \mathscr L_1,,\ldots,\mathscr L_n\right>^1_{(\mathbb X-p_1^{-1}(W))/V}\xrightarrow{\cong}(\pi|_{V})^\ast\left<\mathscr L,\mathscr L_1,\ldots,\mathscr L_n\right>^2_{(\mathbb X-p_1^{-1}(W))/V}\,.
\end{equation}
where $V\subset \mathbb P$ is the image of the flat locus of $p$ (remember that $V$ is open) and $W=\mathbb P-V$. Let's now put $q:=p_2|_{\mathbb X-p_1^{-1}(W)}$. By applying axiom $(1)$, it is enough to get a functorial isomorphism:
\begin{equation}\label{iso_to_show2}
\left<q^\ast\mathscr L,q^\ast\mathscr L_1,\ldots,q^\ast\mathscr L_n\right>^1_{(\mathbb X-p_1^{-1}(W))/V}\xrightarrow{\cong}\left<q^\ast\mathscr L, q^\ast\mathscr L_1,\ldots,q^\ast\mathscr L_n\right>^2_{(\mathbb X-p_1^{-1}(W))/V}\,.
\end{equation}
Now remember that by definition of $\mathscr L_f$ we have:
$$q^\ast\mathscr L=\mathscr L_f|_{(\mathbb X-p_1^{-1}(W))}\otimes (p_1^\ast\mathscr O_\mathbb P(1))^{-1}$$ 
Let's put for simplicity of notations $\mathscr M:=\mathscr L_f|_{(\mathbb X-p_1^{-1}(W))}$; by multi-additivity and axiom $(3)$ we only need to find a functorial isomorphism
\begin{equation}\label{iso_to_show3}
\left<\mathscr M,q^\ast\mathscr L_1\ldots,q^\ast\mathscr L_n\right>^1_{(\mathbb X-p_1^{-1}(W))/V}\xrightarrow{\cong}\left<\mathscr M,q^\ast\mathscr L_1\ldots,q^\ast\mathscr L_n\right>^2_{(\mathbb X-p_1^{-1}(W))/V}\,.
\end{equation}
At this point put $Z'(\sigma):=Z(\sigma)-p_1^{-1}(W)$; thanks to axiom $(2)$, it is enough to find a functorial isomorphism: 
\begin{equation}\label{iso_to_show4}
\left<(q^\ast\mathscr L_1)|_{Z'(\sigma)},\ldots,(q^\ast\mathscr L_n)|_{Z'(\sigma)}\right>^1_{Z'(\sigma)/V}\xrightarrow{\cong}\left<(q^\ast\mathscr L_1)|_{Z'(\sigma)},\ldots,(q^\ast\mathscr L_n)|_{Z'(\sigma)}\right>^2_{Z'(\sigma)/V}\,.
\end{equation}
The relative dimension of the map $p:Z'(\sigma)\to V$ is now $n-1$ and we can apply the inductive hypothesis.

We still have to prove the the existence of $\Psi(\mathscr L_0,\ldots, \mathscr L_n)$ for a general $\mathscr L_0$. For any invertible sheaf $\mathscr L_0$ there exists a sufficiently ample one $\mathscr M$ such that $\mathscr L_0\otimes\mathscr M$ is again sufficiently ample. So we can put:
$$\Psi(\mathscr L_0,\ldots, \mathscr L_n):=\Psi'(\mathscr L_0\otimes \mathscr M,\ldots, \mathscr L_n)\otimes \Psi'(\mathscr M,\ldots, \mathscr L_n)^{-1}$$
provided that the construction doesn't depend on the choice of $\mathscr M$. Such a claim is equivalent to show that $\Psi'(\cdot\,,\mathscr L_1,\ldots, \mathscr L_n)$ is additive with respect to sufficiently ample invertible sheaves.

Now we consider two sufficiently ample invertible sheaves $\mathscr L^{(i)}$ for $i=1,2$ and the associated diagrams:
\begin{equation}\label{diag_univ1}
\begin{tikzcd}[row sep=large, column sep = huge]
\mathbb X^{(i)}:=\mathbb P^{(i)}\times_S X\arrow["p^{(i)}_2"]{r}\arrow["p^{(i)}_1"]{d} & X\arrow["f"]{d}\\
\mathbb P^{(i)}\arrow["\pi^{(i)}"]{r} & S
\end{tikzcd}
\end{equation}
where clearly $\mathbb P^{(i)}:=\mathbb P_S(\mathscr M^{(i)})$ for $\mathscr M^{(i)}:=(f_\ast(\mathscr L^{(i)}))^\vee$. On the other hand if we put $\mathscr L:=\mathscr L^{(1)}\otimes \mathscr L^{(2)}$ and $\mathbb P:=\mathbb P_S((f_\ast\mathscr L^{(1)})^\vee \otimes (f_\ast\mathscr L^{(2)})^\vee)$ we end up with the diagram (\ref{diag_univ1}). There is a natural map 
$$\iota: \mathbb X^{(1)}\times_X \mathbb X^{(2)}\to\mathbb X$$
Let $q_i$ the projections of $\mathbb X^{(1)}\times_X \mathbb X^{(2)}$ on the factors $\mathbb X^{(i)}$, then one can show that:
$$\iota^\ast\mathscr L_f=q_1^\ast\mathscr L^{(1)}_f\otimes q_2^\ast\mathscr L^{(2)}_f\,,$$
$$\iota^\ast \sigma= q_1^\ast\sigma^{(1)}\otimes q_2^\ast\sigma^{(2)}\,.$$
From the properties of the universal extension discussed in \cite[I.2]{Elk} the claim follows.

It remains to show that $\Psi$ transforms $\alpha^{1},\beta^{1},\gamma^{1},\delta^{1}$ to $\alpha^{2},\beta^{2},\gamma^{2},\delta^{2}$. Let's do it for $\alpha^{i}$, the other cases are similar. In particular we have to prove that, given a base change diagram as in axiom (1),  we get a commutative diagram:
\begin{equation}\label{diag_alpha}
\begin{tikzcd}[row sep=large, column sep = huge]
g^\ast\left<\mathscr L_0,\ldots,\mathscr L_n\right>^{1}_{X/S}\arrow["\cong"]{r}\arrow["\cong"]{d} & \left<g'^\ast \mathscr L_0,\ldots,g'^\ast\mathscr L_n\right>^1_{X'/S'}\arrow["\cong"]{d}\\
g^\ast\left<\mathscr L_0,\ldots,\mathscr L_n\right>^{2}_{X/S}\arrow["\cong"]{r} & \left<g'^\ast \mathscr L_0,\ldots,g'^\ast\mathscr L_n\right>^2_{X'/S'}\,.
\end{tikzcd}
\end{equation}
In order to construct (\ref{diag_alpha}) it is enough to proceed similarly  as we did above: we work by induction on $n$. If $n=0$ the claim follows from the property of $\delta^i$ with respect to base change.  For the generic $n$ we can use the universal extension procedure described above and the properties of $\beta^i$ with respect to base change to reduce to $n-1$.
\endproof

\subsection{Deligne pairing in terms of determinant of cohomology}\label{det_del}
In this section we heavily use the proprieties of the determinant functor presented in appendix   \ref{det} in order to give an explicit expression of the Deligne pairing in terms of the determinant of cohomology.

Let $f:X\to S$ be a flat morphism between Noetherian integral schemes. One immediately notices that $\det Rf_\ast $ descends to a map on $K_0(X)$. Now we put:
$$
\left\{\begin{array} {lll}
\left<\mathscr L_0,\ldots, \mathscr L_n\right>_{X/S}:={\det Rf_\ast\left(c_1(\mathscr L^{-1}_0)c_1(\mathscr L^{-1}_1)\ldots c_1(\mathscr L^{-1}_n)\mathscr O_X\right)}\, & \textrm{if } n>0\\
&\\
\left< \mathscr L_0\right>_{X/S}:={\det Rf_\ast\left(c_1(\mathscr L^{-1}_0)\mathscr O_X\right)}^{-1}=N_{X/S}(\mathscr L_0)& \textrm{if } n=0\\
\end{array}\right.
$$
We want to show that this defines the Deligne pairing i.e. that there are some ``canonical''  natural transformations associated to $\left<\;,\ldots, \;\right>_{X/S}$ satisfying all axioms of definition \ref{axiom_deligne}. 
\begin{remark}
When $n=1$, after some simple algebraic manipulations we obtain the expected result:
$$\left<\mathscr L,\mathscr M\right>=\det Rf_\ast(c_1(\mathscr L^{-1}) c_1(\mathscr M^{-1})\mathscr O_X)=\det Rf_\ast(c_1(\mathscr L^{-1})(\mathscr O_X-\mathscr M))=$$
$$=\det Rf_\ast(\mathscr O_X)\otimes\det Rf_\ast(\mathscr L)^{-1}\otimes\det Rf_\ast(\mathscr M)^{-1}\otimes\det Rf_\ast(\mathscr L\otimes\mathscr M) $$
\end{remark}

Like in the case of algebraic varieties, Proposition \ref{propc_1} ensures that $\left<\;,\ldots, \;\right>_{X/S}$ is multi-monoidal and symmetric. Moreover axiom (4) of Definition \ref{axiom_deligne} is trivially satisfied by definition (see  Definition \ref{norm}). So it remains to show that axioms (1)-(3) are satisfied.

\begin{proposition}[Axiom (1) holds]
For any commutative square given by a base change $g:S'\to S$ which is proper, flat and with connected fibres 

$$
\begin{tikzcd}[row sep=large, column sep = huge]
X'=X\times_S S'\arrow["f'"]{r}\arrow["g'"]{d} & S'\arrow["g"]{d}\\
X\arrow["f"]{r} & S
\end{tikzcd}
$$
there is a natural transformation $\alpha_{f,g}$ between multi-monoidal functors $\catname{Pic}(X)^{n+1}\to \catname{Pic}(S')$ such that
 $$\alpha_{f,g}:g^\ast\left<\mathscr L_0,\ldots,\mathscr L_n\right>_{X/S}\xrightarrow{\cong}\left<g'^\ast \mathscr L_0,\ldots,g'^\ast\mathscr L_n\right>_{X'/S'}\,.$$
\end{proposition}
\proof
First of all we have that for any invertible sheaf $\mathscr L$ on $X$ and any coherent sheaf $\mathscr F'$ on $X'$:
$$c_1(\mathscr L)Rg'_\ast\mathscr F'=Rg'_\ast (c_1(g'^\ast\mathscr L)\mathscr F')$$
(see for example the proof of \cite[Lemma B.15]{Kl} for a detailed explanation of the above equality). Therefore
$$c_1(\mathscr L_0)\ldots c_1(\mathscr L_n)Rg'_\ast(\mathscr O_{X'})=c_1(\mathscr L_0)\ldots c_1(\mathscr L_{n-1})Rg'_\ast (c_1(g'^\ast\mathscr L_n)\mathscr O_{X'})=\ldots$$
$$\ldots=Rg'_\ast(c_1(g'^\ast\mathscr L_0)\ldots c_1(g'^\ast\mathscr L_n)\mathscr O_{X'})\,.$$
It means that:

\begin{equation}\label{functor_mess}
g^\ast\det Rf_\ast\left(c_1(\mathscr L_0)\ldots c_1(\mathscr L_n)Rg'_\ast(\mathscr O_{X'})\right)=g^\ast\det Rf_\ast \left(Rg'_\ast(c_1(g'^\ast\mathscr L_0)\ldots c_1(g'^\ast\mathscr L_n)\mathscr O_{X'})\right)\,.
\end{equation}
But thanks to the properties of the morphism $g:S'\to S$ we have that $g'_\ast\mathscr O_{X'}=\mathscr O_X$ (see for example \cite[Exercise 3.11]{Kol}, so it follows that the left hand side of equation (\ref{functor_mess}) is 
$$g^\ast\det Rf_\ast\left(c_1(\mathscr L_0)\ldots c_1(\mathscr L_n)\mathscr O_X\right)\,.$$
 On the right hand side of equation (\ref{functor_mess}) note that we have the composition of the following functors:
\begin{equation}\label{functor_mess1}
g^\ast\circ\det{\!}_{S}\circ Rf_\ast\circ Rg'_\ast\,.
\end{equation}
By the properties of the determinant functor, equation (\ref{functor_mess1}) is naturally isomorphic to
$$
\det{\!}_{S'}\circ Lg^\ast \circ Rf_\ast\circ Rg'_\ast\cong\det{\!}_S'\circ Rf'_\ast=\det Rf'_\ast\,.
$$
In other words we obtained that the right hand side of equation (\ref{functor_mess}) is naturally isomorphic to 
$$\det Rf'_\ast\left(c_1(g'^\ast\mathscr L_0)\ldots c_1(g'^\ast\mathscr L_n)\mathscr O_{X'}\right)\,.$$
\endproof

\begin{remark}
For axioms (2) and (3), we only have to show that the natural transformations exist, since their ``good behaviour'' with respect to base change is ensured by the properties of the determinant of cohomology with respect to base change i.e. equation (\ref{det_bchange}).
\end{remark}

\begin{proposition}[Axiom (2) holds]
When $n>0$ and $D\in\Div(X)$ is an effective relative Cartier divisor, there is a natural transformation $\beta_{f,D}$ between multi-monoidal functors $\catname{Pic}(X)^{n}\to \catname{Pic}(S)$ such that
  $$\beta_{f,D}:\left<\mathscr O_X(D),\ldots,\mathscr L_n\right>_{X/S}\xrightarrow{\cong}\left<\mathscr L_1|_D,\ldots,\mathscr L_n|_D\right>_{D/S}\,.$$
Moreover such a transformation is natural with respect to  base change.
\end{proposition}
\proof
This follows by the simple fact that $c_1(\mathscr O_X(D))\mathscr O_X=\mathscr O_D$ (see Proposition \ref{propc_1}(iii)).
\endproof

\begin{proposition}[Axiom (3) holds]
When $n>0$ there is a natural transformation $\gamma_{f}$ between multi-monoidal functors $\catname{Pic}(S)\times\catname{Pic}(X)^{n}\to \catname{Pic}(S)$
 such that
  $$\gamma_{f}:\left<f^\ast\mathscr L,\mathscr L_1\ldots,\mathscr L_n\right>_{X/S}\xrightarrow{\cong}\mathscr L^{(\mathscr L_1|_{X_s}.\mathscr L_2|_{X_s}.\ldots\mathscr L_n|_{X_s};\mathscr O_{X_s})}$$
where $X_s$ is a generic fibre of $f$. Moreover such a transformation is natural with respect to  base change.
\end{proposition}
\proof
Let's put $\mathscr E:= c_1(\mathscr L_1)\ldots c_1(\mathscr L_n)\mathscr O_X$, then:
$$
\det Rf_\ast (c_1(f^\ast \mathscr L)c_1(\mathscr L_1)\ldots c_1(\mathscr L_n)\mathscr O_X)=\det Rf_\ast(c_1(f^\ast \mathscr L)\mathscr E)=\det Rf_\ast(\mathscr E-(f^\ast \mathscr L)^{-1}\otimes \mathscr E)=
$$
$$
=\det Rf_\ast(\mathscr E)\otimes \det Rf_\ast((f^\ast \mathscr L)^{-1}\otimes \mathscr E)^{-1}=(\ast)
$$
Now, thanks to Proposition \ref{det_chi} the above chain of equalities can be continued in the following way through a canonical isomorphism:
$$
(\ast)\cong\det Rf_\ast(\mathscr E)\otimes \det Rf_\ast(\mathscr E)^{-1}\otimes \mathscr L^{\chi_S(\mathscr E|_{X_s})}=\mathscr L^{\chi_S(\mathscr E|_{X_s})}
$$
where $X_s$ is a generic fibre. In order to conclude, it is enough to notice that $\mathscr E|_{X_s}=c_1(\mathscr L_1|_{X_s})\ldots c_1(\mathscr L_n|_{X_s})\mathscr O_{X_s}$. 

\endproof

\begin{appendices}

\section{Determinant functor and determinant of cohomology}\label{det}

In this section we briefly discuss, without proofs, the determinant functor by following \cite{KM}. First we define the determinant for locally free sheaves, then we extend it to complexes of locally free sheaves and then we extend it further for perfect complexes. We will use some basic notions from the theory of derived category (see \cite{Hoche} for a concise introduction.)

We fix a Noetherian integral scheme $X$. A \emph{graded invertible sheaf} on $X$ is a couple $(\mathscr L, \alpha)$ where $\mathscr L$ is an invertible sheaf on $X$ and $\alpha: X\to\mathbb Z$ is a continuous function. A morphism of graded invertible sheaves $\phi: (\mathscr L, \alpha)\to(\mathscr M,\beta)$ is a morphism of invertible sheaves such that the following condition hold: for any $x\in X$, if $\alpha(x)\neq \beta(x)$, then $\phi_x=0$. We denote with $\catname{Gr}(X)$ the category of graded invertible sheaves, and $\catname{isGr}(X)$ is the category whose objects are graded invertible sheaves and the morphisms are just the isomorphisms; note that $\catname{isGr}(X)$ is a Picard groupoid. The tensor product (i.e. the group operation) between graded invertible sheaves is defined as $(\mathscr L, \alpha)\otimes(\mathscr M, \beta)=(\mathscr L\otimes\mathscr M,\alpha+\beta)$. The unit graded invertible sheaf is $(\mathscr O_X,0)$. Furthermore we can define the isomorphism $\tau:\mathscr L\otimes\mathscr M\to\mathscr M\otimes\mathscr L$ such that locally and on pure tensors is given by:
$$\tau(l\otimes m)=(-1)^{\alpha\beta}m\otimes l\,.$$
Let $\catname{Vec}(X)$ be the category of locally free sheaves on $X$ (of finite rank) and let $\catname{isVec}(X)$ its subcategory where the morphisms are only the isomorphisms.  

 For a locally free  sheaf $\mathscr E$ of rank $r$, we  denote  with the symbol $\bigwedge^r\mathscr E$ the sheafification of the following presheaf:
$$U\mapsto\bigwedge^r\mathscr E(U)\,.$$
Then we can construct graded invertible sheaf $\det{\!}^\star\mathscr E$ in the following way:
$$\det{\!}^\star\mathscr E:=\left(\bigwedge^r\mathscr E, r\right)\,.$$
Then we have a functor:
$$\det{\!}_X^\star\colon \catname{isVec}(X)\to\catname{isGr}(X)\,.$$
For any short exact sequence of locally free sheaves:
\begin{equation*}
\begin{tikzcd}
0\arrow[r]&\mathscr H\arrow[r,"\alpha"]&\mathscr F\arrow[r,"\beta"]&\mathscr G\arrow[r]& 0\
\end{tikzcd}
\end{equation*} 
there is an isomorphism of graded invertible sheaves
\begin{equation*}
\begin{aligned}
i_X^\star(\alpha,\beta):\det{\!}_X^\star(\mathscr H)\otimes\det{\!}_X^\star(\mathscr G)\xrightarrow{\cong} \det{\!}_X^\star(\mathscr F)
\end{aligned}
\end{equation*}
that locally is given in the following way: assume that $\mathscr H$ has rank $r$ and $\mathscr G$ has rank $s$, then for any local sections $h_i\in\mathscr H(U)$ and $\beta(f_i)\in\mathscr G(U)$, for $f_i\in\mathscr F(U)$ we have:

$$i_X^\star(\alpha,\beta)((h_1\wedge\ldots\wedge h_r)\otimes(\beta(f_1)\wedge\ldots\wedge\beta(f_s)))=\alpha(h_1)\wedge\ldots\wedge \alpha(h_r)\wedge f_1\wedge\ldots\wedge f_n\,.$$
We are ready to give the definition of determinant for bounded complexes of locally free sheaves.
\begin{definition}\label{def_det}
Let $\catname{isVec}_b^\bullet(X)$ be the category of bounded complexes in $\catname{Vec}(X)$ where the morphism are just the quasi-isomorphisms between complexes. Then a \emph{determinant functor on $X$} consists of the following data:
\begin{enumerate}
\item[$(1)$] A functor $\mathfrak f_X:\catname{isVec}_b^\bullet(X)\to\catname{isGr(X)}$
\item[$(2)$] For any short exact sequence in $\catname{isVec}_b^\bullet(X)$:
\begin{equation*}
\begin{tikzcd}
0\arrow[r]&\mathscr H^\bullet\arrow[r,"\alpha"]&\mathscr F^\bullet\arrow[r,"\beta"]&\mathscr G^\bullet\arrow[r]& 0\
\end{tikzcd}
\end{equation*}
an isomorphism:

\begin{equation*}
\begin{aligned}
i_X(\alpha,\beta):\mathfrak f_X(\mathscr H^\bullet)\otimes\mathfrak f_X(\mathscr G^\bullet)\xrightarrow{\cong} \mathfrak f_X(\mathscr F^\bullet)
\end{aligned}
\end{equation*}

\end{enumerate} 

Moreover $(\mathfrak f_X, i_X)$ satisfies the following conditions:
\begin{enumerate}
\item[$(i)$] Given a commutative diagram in $\catname{isVec}_b^\bullet(X)$:
\begin{equation*}
\begin{tikzcd}
0\arrow[r]&\mathscr H^\bullet\arrow[r, "\alpha"]\arrow[d,"\lambda'"]&\mathscr F^\bullet\arrow[r, "\beta"]\arrow[d, "\lambda"]&\mathscr G^\bullet\arrow[r]\arrow[d,"\lambda''"]& 0\\
0\arrow[r]&\mathscr H_1^\bullet\arrow[r, "\alpha_1"]&\mathscr F_1^\bullet\arrow[r, "\beta_1"]&\mathscr G_1^\bullet\arrow[r]& 0\
\end{tikzcd}
\end{equation*}
such that the rows are short exact sequences, then the following diagram commutes
$$
\begin{tikzcd}[row sep=large, column sep = large]
\mathfrak f_X({\mathscr H}^\bullet)\otimes \mathfrak f_X({\mathscr G}^\bullet)\arrow["{i_X(\alpha,\beta)}"]{r}\arrow["\mathfrak f_X(\lambda')\otimes \mathfrak f_X(\lambda'')"]{d} & \mathfrak f_X({\mathscr F}^\bullet)\arrow["\mathfrak f_X(\lambda)"]{d}\\
\mathfrak f_X({\mathscr H_1}^\bullet)\otimes \mathfrak f_X({\mathscr G_1}^\bullet)\arrow["{i_X(\alpha_1,\beta_1)}"]{r} & \mathfrak f_X({\mathscr F_1}^\bullet)
\end{tikzcd}
$$

\item[$(ii)$] Given a commutative diagram in $\catname{isVec}_b^\bullet(X)$:
$$
\begin{tikzcd}
  & 0\arrow{d} & 0\arrow{d} & 0\arrow{d} &\\
  0 \arrow{r} & \mathscr H_1^\bullet \arrow[r,"\alpha_1"]\arrow[d,"\gamma_1"] & \mathscr F_1^\bullet \arrow[r,"\beta_1"]\arrow[d,"\gamma"] & \mathscr G_1^\bullet \arrow{r}\arrow[d,"\gamma_2"] & 0\\
  0 \arrow{r} & {\mathscr H}^\bullet \arrow[r,"\alpha"]\arrow[d,"\delta_1"] & {\mathscr F}^\bullet \arrow[r,"\beta"]\arrow[d, "\delta"] & {\mathscr G}^\bullet \arrow{r}\arrow[d,"\delta_2"] & 0
\\
 0 \arrow{r} & {\mathscr H_2}^\bullet \arrow[r,"\alpha_2"]\arrow[d] & {\mathscr F_2}^\bullet \arrow[r,"\beta_2"]\arrow[d] & {\mathscr G_2}^\bullet \arrow{r}\arrow[d] & 0\\
 & 0 & 0 & 0 &
\end{tikzcd}
$$
such that  rows and columns are short exact sequences,  then the following diagram commutes
$$
\begin{tikzcd}[row sep=huge, column sep = 4.5cm]
\mathfrak f_X({\mathscr H_1}^\bullet)\otimes \mathfrak f_X({\mathscr H_2}^\bullet)\otimes\mathfrak f_X({\mathscr G_1}^\bullet)\otimes \mathfrak f_X({\mathscr G_2}^\bullet)\arrow["{i_X(\gamma_1,\gamma_2)\otimes i_X(\delta_1,\delta_2)}"]{r}\arrow["{i_X(\alpha_1,\beta_1)\otimes i_X(\alpha_2,\beta_2)}"]{d} & \mathfrak f_X({\mathscr H}^\bullet)\otimes \mathfrak f_X({\mathscr G}^\bullet)\arrow["{i_X(\alpha,\beta)}"]{d}\\
\mathfrak f_X({\mathscr F_1}^\bullet)\otimes \mathfrak f_X({\mathscr F_2}^\bullet)\arrow["{i_X(\gamma,\delta)}"]{r} & \mathfrak f_X({\mathscr F}^\bullet)
\end{tikzcd}
$$
\item[$(iii)$] $\mathfrak f_X$ and $i_X(\,,\,)$ both commute with the base change of $X$. The explicit expression of such a property is the following: fix  a morphism of schemes $\psi:Y\to X$ and let $L\psi^\ast:\catname D_-(\catname{ QCoh}(X))\to\catname D_-(\catname {QCoh}(Y))$ be the left derived functor of to the pullback $\psi^\ast$; then the following properties hold:
\begin{enumerate}
\item[\sbt] There is a natural transformation between functors: $\eta(\psi)\colon\mathfrak f_Y\circ L\psi^\ast\xrightarrow{\cong}\psi^\ast\circ\mathfrak f_X$.
\item[\sbt] For any short exact sequence in $\catname{isVec}_b^\bullet(Y)$:
\begin{equation*}
\begin{tikzcd}
0\arrow[r]&\mathscr H^\bullet\arrow[r,"\alpha"]&\mathscr F^\bullet\arrow[r,"\beta"]&\mathscr G^\bullet\arrow[r]& 0\
\end{tikzcd}
\end{equation*}
then the following diagram commutes
$$
\begin{tikzcd}[row sep=large, column sep = huge]
\mathfrak f_Y(L\psi^\ast{\mathscr H}^\bullet)\otimes\mathfrak f_Y(L\psi^\ast{\mathscr G}^\bullet)\arrow["{i_Y(L\psi^\ast(\alpha,\beta))}"]{r}\arrow["\eta(\psi)(\mathscr H^\bullet)\otimes \eta(\psi)(\mathscr G^\bullet)"]{d} & \mathfrak f_Y(L\psi^\ast{\mathscr F}^\bullet)\arrow["\eta(\psi)(\mathscr F^\bullet)"]{d}\\
\psi^\ast\mathfrak f_X({\mathscr H}^\bullet)\otimes \psi^\ast\mathfrak f_X({\mathscr G}^\bullet)\arrow["{\psi^\ast i_X(\alpha,\beta)}"]{r} & \psi^\ast\mathfrak f_X({\mathscr F}^\bullet)
\end{tikzcd}
$$
\end{enumerate}

\item[$(iv)$] $\mathfrak f_X(0)=(\mathscr O_X, 0)$. Moreover for the short exact sequence:
\begin{equation*}
\begin{tikzcd}
0\arrow[r]&\mathscr F^\bullet\arrow[r,"\id"]&\mathscr F^\bullet\arrow[r,"0"]&0\arrow[r]& 0\
\end{tikzcd}
\end{equation*}
we have that $i_X(\id,0)\colon \mathfrak f_X(\mathscr F^\bullet)\otimes(\mathscr O_X,0)\to  \mathfrak f_X(\mathscr F^\bullet)$ is the canonical map (i.e. the ``projection on the first component'').
\item[$(v)$] If we canonically identify $\catname{isVec}(X)$ as a subcategory of $\catname{isVec}_b^\bullet(X)$, then $\mathfrak f_X$ restricts to $\det{\!}_X^\star$ and $i_X(\,,\,)$ restricts to $i^\star_X(\,,\,)$.
\end{enumerate}
\end{definition}

\begin{theorem}
Up to natural transformation of functors there exists a unique determinant $(\mathfrak f_X, i_X)$ on $X$. It will be denoted as $(\det_X, i_X)$.
\end{theorem}
\proof
See \cite[Theorem 1]{KM} for a complete proof. Here we just write down the explicit expressions for $\det_X$:
\begin{equation}
\det{\!}_X(\mathscr F^\bullet)=\bigotimes_i\det{\!}^\star_X(\mathscr F^i)^{(-1)^i}\,.
\end{equation}

\endproof

The category   $\catname{isVec}_b^\bullet(\cdot)$ is quite restrictive, for example it doesn't behave well with respect to the pushforward functor. Therefore,  we would like to have a determinant functor for a more general category.
\begin{definition}
A complex $\mathscr F^\bullet $ of $\mathscr O_X$-modules is said \emph{perfect} if for any $x\in X$ there exist an open neighbourhood $U\ni x$, a complex $\mathscr G^\bullet $ in $\catname {Vec}^\bullet_b(U)$ and a quasi isomorphism  of complexes of $\mathscr O_U$-modules $\mathscr G^\bullet\to\mathscr F^\bullet_{|U}$. The category of perfect complexes on $X$ is denoted by $\catname{Perf}(X)$, whereas $\catname{isPerf}(X)$ denotes the category of perfect complexes where the morphisms are just the quasi-isomorphisms.
\end{definition}
\begin{theorem}
The determinant $(\det_X,i_X)$ can be extended uniquely, up to natural transformation, to a determinant on $\catname{isPerf}(X)$. We will denote this extension again with the symbol $(\det_X,i_X)$ and formally it is the datum of:
\begin{enumerate}
\item[$(1)$] A functor $\det_X\colon\catname{isPerf}(X)\to\catname{isGr(X)}$
\item[$(2)$] For any short exact sequence in $\catname{isPerf}(X)$:
\begin{equation*}
\begin{tikzcd}
0\arrow[r]&\mathscr H^\bullet\arrow[r,"\alpha"]&\mathscr F^\bullet\arrow[r,"\beta"]&\mathscr G^\bullet\arrow[r]& 0\
\end{tikzcd}
\end{equation*}
an isomorphism:
\begin{equation*}
\begin{aligned}
i_X(\alpha,\beta):\det{\!}_X(\mathscr H^\bullet)\otimes\det{\!}_X(\mathscr G^\bullet)\xrightarrow{\cong} \det{\!}_X (\mathscr F^\bullet)
\end{aligned}
\end{equation*}
\end{enumerate} 
Moreover the properties $(i)-(v)$ listed in definition \ref{def_det} are satisfied in $\catname{isPerf}(X)$. 
\end{theorem}
\proof
See \cite[Theorem 2]{KM}.
\endproof
One of the most important applications of the determinant functor appears in arithmetic geometry if we consider its interaction with the usual pushfoward functor.

\begin{definition}
Let $f:X\to S$ be a flat morphism between integral Noetherian schemes and let $\mathscr F$ be a coherent sheaf on $X$. It is well known (\cite[Exp. 3, Proposition 4.8]{SGA}) that the complex $Rf_\ast\mathscr F$ induced by the right derived functor of $f_\ast$ is a perfect complex on $S$. Then, just by composing $Rf_\ast$ with $\det_S$ it is possible to define the functor: 
$$\det Rf_\ast:=\det{\!}_S\circ Rf_\ast : \catname{isCoh}(X)\to \catname{isGr}(S)$$
which is called the \emph{determinant of cohomology (relative to $f$)}. Very often, for  simplicity  we want to forget about the graduation on the target of the determinant of cohomology, so it becomes a functor  $\catname{isCoh}(X)\to \catname{Pic}(S)$.
\end{definition}

Since the right derived functor $Rf_\ast$ is exact in the derived sense (see \cite{Hoche}) it is not hard to show that for any short exact sequence of coherent sheaves on $X$
\begin{equation*}
\begin{tikzcd}
0\arrow[r]&\mathscr H\arrow[r,"\alpha"]&\mathscr F\arrow[r,"\beta"]&\mathscr G\arrow[r]& 0\
\end{tikzcd}
\end{equation*} 
there is an isomorphism of graded invertible sheaves
\begin{equation*}
\begin{aligned}
i_f(\alpha,\beta):\det Rf_\ast\mathscr H\otimes\det Rf_\ast\mathscr G\xrightarrow{\cong} \det Rf_\ast\mathscr F
\end{aligned}
\end{equation*}
Moreover the whole construction behaves well with respect to flat base change in the following sense: assume that the following commutative square is given by a flat base change from $S$ to $S'$

$$
\begin{tikzcd}[row sep=large, column sep = huge]
X'=X\times_S S'\arrow["f'"]{r}\arrow["g'"]{d} & S'\arrow["g"]{d}\\
X\arrow["f"]{r} & S
\end{tikzcd}
$$
then for any coherent sheaf $\mathscr F$ on $X$ we have

\begin{equation}\label{det_bchange}
g^\ast(\det Rf_\ast \mathscr F)\cong \det Rf'_\ast(g'^\ast\mathscr F)\,.
\end{equation}
We will also need an important property of the determinant of cohomology:

\begin{proposition}\label{det_chi}
Let $f:X\to S$ be a flat morphism between integral Noetherian schemes and let $\mathscr F$ be a coherent sheaf on $X$. Moreover let $\mathscr L$ be an invertible sheaf on $S$. Then there is a canonical isomorphism between invertible sheaves on $S$:
$$\det Rf_\ast(f^\ast \mathscr L\otimes \mathscr F)\xrightarrow{\cong} \mathscr L^{\otimes \chi_S(\mathscr F|_{X_s})}\otimes \det Rf_\ast(\mathscr F)$$
where $X_s$ is a generic fibre of $f$ .
\end{proposition}
\proof
See \cite{KM}.
\endproof

\begin{definition}\label{norm}
Let $\varphi:X\to S$ be a finite morphism between integral Noetherian schemes, then \emph{the norm of $\varphi$} is defined as the functor:
\begin{eqnarray*}
N_{\varphi}=N_{X/S}:\catname{Pic}(X) &\to & \catname{Pic}(S)\\
\mathscr L &\mapsto & \det R\varphi_\ast \mathscr L\otimes (\det R\varphi_\ast \mathscr O_X)^{-1}\,.
\end{eqnarray*} 
\end{definition}

\section{Original construction of Deligne pairing}\label{ori_constr}
In this section we give all details of the construction of the Deligne pairing described in \cite{Del}.

$S$ is a Dedekind scheme\footnote{For us a Dedekind scheme is an integral, Noetherian, normal, scheme of dimension 0 or 1.} and we put $K:=K(S)$.  $\varphi:X\to S$ is a $S$-scheme satisfying the following properties:
\begin{enumerate}
\item[\sbt] $X$ is  two dimensional, integral, and regular. The generic point of $X$ is $\eta$ and the function field of $X$ is denoted by $K(X)$.
\item[\sbt] $\varphi$ is proper and flat.
\item[\sbt] The generic fibre, denoted by $X_K$, is a geometrically integral, smooth, projective curve over $K$. 
\end{enumerate} 
We say that $X$ is an \emph{arithmetic surface over $S$}. 

In this section we will also need recall the norm operator $\mathcal N$ in dimension $1$ and $2$ (it is formally different from the norm of an invertible sheaf defined above). 
\begin{definition}
let $C$ be a projective, non-singular curve over a field $k$, then for a closed point $x\in C$  and any non-zero rational function $f\in K(C)^\times$ such that $f\in\mathscr O^\times_{C,x}$ we put
\begin{equation}\label{N_oncurves}
\mathcal N_x(f):=N_{k(x)|k} \left(f(x)\right)\,,
\end{equation}
where $f(x)$ is the obvious element of $k(x)$ associated to $f$. So if $D=\sum_{x\in C}n_x[x]\in \Div(C)$ and $f\in K(C)^\times$ is a non-zero rational function such that $(f)$ and $D$ have no common components, then it is well defined the following element:
$$ \mathcal N_D(f):=\prod_{x\in C}  \mathcal N_x(f)^{n_x}\quad\in k^\times$$
\end{definition}
The well known Weil reciprocity law  says that:
$$\mathcal N_{(g)}(f)=\mathcal N_{(f)}(g)\rlap.$$

\noindent Coming back to our arithmetic surface $\varphi:X\to S$, consider 
 $$\Upsilon:=\Set{(D,E)\in\Div(X)\times\Div(X)\colon\textrm{$D$ and $E$ have no common components}}\rlap,$$
and note that if $(D_j,E_j)\in\Upsilon$ with $j=1,2$, then $(D_1+D_2,E_1+E_2)\in\Upsilon$. 
\begin{definition}\label{def1}
Let $(D,E)\in \Upsilon$ such that $D$ and $E$ are both effective, then for any closed point $x\in X$ we put:

$$i_x(D,E):= \len_{\mathscr O_{X,x}}\mathscr O_{X,x}/\left (\mathscr O_X(-D)_x+\mathscr O_X(-E)_x\right)\,.$$
This is called the \emph{local intersection number} of $D$ and $E$ at $x$. 
\end{definition}
The local intersection number assigns the multiplicity of the intersection at each point of $X$,  and the following basic result summarizes its naive properties. 
\begin{proposition}\label{prop_loc_int}
Let $(E,D)\in \Upsilon$ and $(E_j,D_j)\in\Upsilon$ with $j=1,2$ such that all the divisors are effective, then
\begin{enumerate} 
\item[$(1)$]$i_x(D,E)=i_x(E,D)$.
\item[$(2)$]$\displaystyle i_x(D_1+D_2,E_1+E_2)=\sum_{j,k=1}^2 i_x(D_j,E_k)$.
\item[$(3)$] $i_x(D,E)\neq 0$ if and only if $x\in\supp(D)\cap\supp(E)$.
\item[$(4)$] If $x\in E$, $i_x(D,E)=\mult_x(D|_E)$.
\end{enumerate}
\end{proposition}
\proof
$(1)$ and $(3)$ are obvious. For $(2)$ and $(4)$ see \cite[lemma 9.1.4]{Liu}.
\endproof
Any divisor $D\in\Div(X)$ can be written in a unique way as $D=D_+-D_-$ where both $D_+$ and $D_-$ are effective and  if $(D,E)\in\Upsilon$, then $(D_\pm,E_\pm)\in\Upsilon$. We can use definition \ref{def1} in order to  have the local intersection at $x$ of $D$ and $E$ when $(D,E)$ is any element of $\Upsilon$ (so not necessarily effective):
$$i_x(D,E):=i_x(D_+,E_+)-i_x(D_+,E_-)-i_x(D_-,E_+)+i_x(D_-,E_-)\,.$$
\begin{definition}\label{def2}
Let $(D,E)$ be an element of $\Upsilon$, then we define the $0$-cycle on $X$ given by:
$$i(D,E):=\sum_{x\in X^{(0)}}i_x(D,E)[x]\rlap,$$
where here $[x]$ is a shorthand of $[\overline{\{x\}}].$
\end{definition}
\begin{remark}
The sum in  definition \ref{def2} is finite because if $D$ and $E$ are effective without common components, then $i_x(D,E)=\mult_x(D|_E)$ (proposition \ref{prop_loc_int}(4)) and there is only a finite number of points on $E$ at which the divisor $D|_E$ has non-zero multiplicity. 
\end{remark}
\begin{proposition}
If $(D,E),(D_j,E_j)\in\Upsilon$ with $j=1,2$, then the following properties hold for $i(D,E)$:
\begin{enumerate}
\item[\sbt] $i(D,E)=i(E,D)$ (symmetry)\,.
\item[\sbt] $\displaystyle i(D_1+D_2,E_1+E_2)=\sum_{j,k=1}^2 i(D_j,E_k)$ (bilinearity)\,.
\end{enumerate}
\end{proposition}
\proof
It follows immediately from proposition \ref{prop_loc_int}.
\endproof
\begin{definition}\label{finite_int_numb}
We have the symmetric and bilinear pairing on $\Upsilon$:
\begin{eqnarray*}
\Upsilon &\to& \Div(S)\\
(D,E) &\mapsto& \left<D,E\right>
\end{eqnarray*}
where
$$\left<D,E\right>:=\varphi_\ast i(D,E)=\sum_{x\in X} [k(x):k(\varphi(x))]\,i_{x}(D,E)\, [\varphi(x)]\,.$$
\end{definition}

Let $\Gamma$ be a prime divisor of $X$ with generic point $\gamma$ and consider a non-zero rational function $f\in K(X)^\times$ such that $(f)$ and $\Gamma$ have no common components, then define $\mathcal N_{\Gamma}(f)\in K^\times$ in the following way:
$$\mathcal N_{\Gamma}(f):=\left\{\begin{array} {cc}
N_{K(\Gamma)|K}(f|_{\Gamma}) & \textrm{if $\Gamma$ is horizontal}\\
1 & \textrm{if $\Gamma$ is vertical}\\
\end{array}\right.
$$
where $N_{K(\Gamma)|K}$ is the usual field norm and $f|_{\Gamma}$ is defined as follows: since $(f)$ and $\Gamma$ have no common components it follows that $v_{\gamma}(f)=0$, that is $f\in\mathscr O_{X,\gamma}^\times$. So $f|_{\Gamma}$, is the natural image of $f$ in $k(\gamma)=K(\Gamma)$. At this point for any  $D=\sum_i n_i\Gamma_i\in \Div(X)$ such that $D$ and $(f)$ have no common components we have:
$$\mathcal N_{D}(f):=\prod_i\mathcal N_{\Gamma_i}(f)^{n_i}\quad \in K^\times$$
Since $K(X)$ is the function field of any open subscheme $U\subseteq X$ and of $X_K$ we can restrict the operator $\mathcal N_\ast(\cdot)$ to $U$ and to $X_K$.

\begin{proposition}\label{restr_prop}
Let $f\in K(X)^\times$ and let $D\in Div(X)$ such that $(f)$ and $D$ have no common components, then the following claims hold:
\begin{itemize}
\item[$(1)$] Let $U\subseteq X$ be an open subscheme, then $\mathcal N_{D|_U}(f)=\mathcal N_D(f)$.

\item[$(2)$] $\mathcal N_{D|_{X_K}}(f)=\mathcal N_D(f)$, where the left hand side is the one-dimensional operator defined in equation (\ref{N_oncurves}).

\end{itemize}
\end{proposition}
\proof
In both items we can restrict to the case when $D=\Gamma$ is an irreducible horizontal divisor.\\
$(1)$ The function fields  and the generic points of $\Gamma$ and $\Gamma|_U$ coincide, so the claim follows trivially.\\
$(2)$ Let $\gamma\in X_K$ be the generic point of $\Gamma$, it is a closed point of $X_K$ such that $k(\gamma)=K(\Gamma)$. By the bare definitions we can check the required equality. 
\endproof

\begin{proposition}\label{intwithrat}
Let $f\in K(X)^\times$ and let $D\in\Div(X)$ a divisor such that $D$ and $(f)$ have no common components, then
$$\left<D,(f)\right>=\left(\mathcal N_D(f)\right)\;\in \Princ(S)\,.$$ 
\end{proposition}
\proof
See \cite[Proposition 4.3]{Mor}.
\endproof

Now we will construct the Deligne pairing and see the relation with the pairing $\left<D,E\right>$ for divisors. We divide the construction in two steps:

\emph{Step 1.} Definition of the $K$-vector space $\left<\mathscr L,\mathscr M\right>_K$.\\

Consider the sets:

\begin{alignat*}{2}
  \Upsilon_K &:= \biggl\{(D,E)\in \Div(X)\times\Div(X)   &&\;\colon\; \pctext{2in}{$D|_{X_K}$ and  $E|_{X_K}$ have no common components (as divisors on $X_K$)}\biggr\}\rlap, 
\end{alignat*}
\begin{alignat*}{2}
  \Sigma_K &:= \biggl\{(l,m)   &&\;\colon\; \pctext{3.5in}{ $l$ and $m$ are non-zero meromorphic sections of $\mathscr L$ and $\mathscr M$ such that $(\divi(l),\divi(m))\in \Upsilon_K$}\biggr\}\rlap. 
\end{alignat*}
Note that $\Upsilon_K$ is  just the set of couple of divisors with no common horizontal components. Now we define some vector spaces over $K$.
$$V:= K^{(\Sigma_K)}\,,$$
namely $V$ is the free $K$-vector space over $\Sigma_K$. 
\begin{equation}\label{rel1}
W':=\left\{(fl,m)-\mathcal N_{\divi(m)|_{X_K}}(f)\cdot(l,m)\colon f\in K(X)^\times, (l,m),(fl,m)\in \Sigma_K\right\}\,,
\end{equation} 

\begin{equation}\label{rel2}
T':=\left\{(l,gm)-\mathcal N_{\divi(l)|_{X_K}}(g)\cdot(l,m)\colon g\in K(X)^\times,  (l,m),(l,gm)\in \Sigma_K \right\}\,.
\end{equation} 
Note that the above ``$\mathcal N_{\ast}(\cdot)$'' is the one-dimensional operator of definition \ref{N_oncurves} considered on the curve $X_K$.

\begin{remark}
$\mathcal N_{\divi(m)|_{X_K}}(f)$ and $\mathcal N_{\divi(l)|_{X_K}}(g)$ are well defined since $(l,m)$, $(fl,m)$, $(l,gm) \in \Sigma_K$, so $\divi(m)|_{X_K}$ and $(f)$ have no common components. The same holds for  $\divi(l)|_{X_K}$ and $(g)$.
\end{remark}

Define the free vector spaces $W:= K^{(W')}$ and $T:=K^{(T')} $; moreover put 
$$\left<\mathscr L,\mathscr M\right>_K:=V/(W+T)\,\rlap,$$
which is considered as a constant sheaf (of $K$-vector spaces) over $X$. The natural image of any element $(l,m)\in \Sigma_K\subset V$ in $\left <\mathscr L,\mathscr M\right>_K$ is denoted  as $\left <l,m\right >_K$.\\
\begin{proposition}\label{sstep1}
 $\left<\mathscr L,\mathscr M\right>_K$ is a one-dimensional vector space over $K$.
\end{proposition}
\proof
Fix $(l_0,m_0)\in\Sigma_K$, then for any $(l,m)\in\Sigma_K$ there are two elements $f_0,g_0\in K(X)^\times$ such that $l=f_0l_0$, $m=g_0m_0$ and moreover:
$$((f_0),(g_0)),\;((f_0),\;\divi(m_0)),\;((g_0),\divi(l_0))\in\Upsilon_K\,.$$
By equations (\ref{rel1}) and (\ref{rel2}), in $\left<\mathscr L,\mathscr M\right>_K$ we can write:
\begin{equation}\label{rel3}
\left<l,m\right>_K= \left<f_0l_0,g_0m_0\right>_K=[f_0,g_0]\,\mathcal N_{\divi(m_0)|_{X_K}}(f_0)\,\mathcal N_{\divi(l_0)|_{X_K}}(g_0)\left<l_0,m_0\right>_K\rlap.
\end{equation}
where, in order  to simplify the notations, we put $[f_0,g_0]:=\mathcal N_{(f_0)}(g_0)$  intended as operation on the curve $X_K$. This shows that $\left<\mathscr L,\mathscr M\right>_K$ has dimension at most $1$ over $K$.
Define the  homomorphism of $K$-vector spaces:
$$\theta:V\to K$$
such that
$$\theta(l,m):=[f_0,g_0] \mathcal N_{\divi(m_0)|_{X_K}}(f_0) N_{\divi(l_0)|_{X_K}}(g_0)\,.$$
 Note that  $\theta$ is non-trivial, so surjective, since $\theta(l_0,m_0)=1$. Now by using the Weil reciprocity law  we prove that $\theta$ descends to a non-trivial morphism $\overline\theta:\left<\mathscr L,\mathscr M\right>_K\to K$, indeed for $f,g\in K(X)^\times$:
\begin{eqnarray*}
\theta(fl,m)&=&[ff_0,g_0]\, \mathcal N_{\divi(m_0)|_{X_K}}(ff_0)\, \mathcal N_{\divi(l_0)|_{X_K}}(g_0)=\\
&=& [f,g_0]\,[f_0,g_0]\, \mathcal N_{\divi(m_0)|_{X_K}}(f)\, \mathcal N_{\divi(m_0)|_{X_K}}(f_0)\, \mathcal N_{\divi(l_0)|_{X_K}}(g_0)=\\
&=&[g_0,f]\, \mathcal N_{\divi(m_0)|_{X_K}}(f)\, \theta(l,m)=\\
&=& \mathcal N_{\divi(m)|_{X_K}}(f)\, \theta(l,m)\rlap.
\end{eqnarray*}
Similarly it holds that
$$\theta(l,gm)=\mathcal N_{\divi(l)|_{X_K}}(g)\, \theta(l,m)\rlap.$$
In other words equation \ref{rel3} can we written as:
$$
\left<l,m\right>_K=\overline{\theta}(\left<l,m\right>_K)\left<l_0,m_0\right>_K
$$
hence, by the  non triviality of $\overline{\theta}$ we conclude that $\left<\mathscr L,\mathscr M\right>_K$ has dimension $1$.
\endproof
\emph{Step 2.} Definition of $\left<\mathscr L,\mathscr M\right>$.\\
 Let $U\subseteq S$ be a non-empty open subset and denote  with $X_U$ the schematic inverse image of $U$ with respect to $\varphi$. We  clearly have a flat map $X_U\to U$, so we define:
\begin{alignat*}{2}
  \Upsilon_U &:= \biggl\{(D,E)\in \Div(X)\times\Div(X)   &&\;\colon\; \pctext{2in}{$D|_{X_U}$ and  $E|_{X_U}$ have no common components (as divisors on $X_U$)}\biggr\}\rlap, 
\end{alignat*}
\begin{alignat*}{2}
  \Sigma_U &:=\Bigg\{(l,m)   &&\;\colon\; \pctext{3in}{$l$ and $m$ are non-zero meromorphic sections of $\mathscr L$ and $\mathscr M$ such that  $(\divi(l),\divi(m))\in\Upsilon_U$  and $\left<\divi(l)|_{X_U}, \divi(m)|_{X_U} \right>$ is effective on $U$}\Bigg\}\rlap. 
\end{alignat*}
Moreover notice that  $\Sigma_U\subset \Sigma_K$. We define  a sheaf of $\mathscr O_S$-modules $\mathscr A$ on $X$ given by:
$$\mathscr A|_{U}:={\mathscr O_{S}|_{U}}^{(\Sigma_{U})}\,.$$
Finally consider the morphism of sheaves: $\Phi:\mathscr A\to\left<\mathscr L,\mathscr M\right>_K$ which sends $(l,m)\in\Sigma_U$ to $\left<l,m\right>_K$ and define
$$\left<\mathscr L,\mathscr M \right>:=\faktor{\mathscr A}{\ker(\Phi)}\rlap.$$
The canonical image of $(l,m)\in\Sigma_U$ in $\mathscr A(U)$ is denoted as $\left<l,m \right>_U$.

\begin{proposition}\label{sstep2}
Let $(l,m)\in \Sigma_{U}$ such that $\left<\divi(l)|_{X_{U}}, \divi(m)|_{X_{U}}\right>=0\in \Div(U)$. Then for any $(l',m')\in\Sigma_{U}$ there exists an element $a\in\mathscr O_S(U)$ such that $\left<l',m'\right>_{U}=a\left<l,m\right>_{U}$.
\end{proposition}
\proof

There are two elements $f,g\in K(X)^\times$ such that $l'=fl$, $m'=gm$ and moreover:
$$((f),(g)),\;((f),\;\divi(m)),\;((g),\divi(l))\in\Upsilon_K\,.$$
Hence by using proposition \ref{intwithrat}:
$$\left<\divi(l')|_{X_{U}}, \divi(m')|_{X_{U}}\right>=\left<(f)|_{X_{U}}\divi(l)|_{X_{U}}, (g)|_{X_{U}}\divi(m)|_{X_{U}}\right>=$$
$$=\left< (f)|_{X_{U}},(g)|_{X_{U}}\right>+\left< (f)|_{X_{U}},\divi(m)|_{X_{U}}\right>+\left<\divi(l)|_{X_{U}},  (g)|_{X_{U_i}}\right>+0=$$
$$= \left(\mathcal N_{(f)|_{X_{U}}}(g)\right)+\left(\mathcal N_{\divi(m)|_{X_{U}}}(f)\right)+\left(\mathcal N_{\divi(l)|_{X_{U}}}(g)\right)=$$
$$=\left(\mathcal N_{(f)|_{X_{U}}}(g)\, \mathcal N_{\divi(m)|_{X_{U}}}(f)\, \mathcal N_{\divi(l)|_{X_{U}}}(g) \right)\rlap.$$
Since $\left<\divi(l')|_{X_{U}}, \divi(m')|_{X_{U}}\right>$ is effective, then 
$$a:=\mathcal N_{(f)|_{X_{U}}}(g)\,\mathcal N_{\divi(m)|_{X_{U}}}(f)\, \mathcal N_{\divi(l)|_{X_{U}}}(g)\in \mathscr O_S(U)\,.$$
On the other hand
$$\left<l',m'\right>_K=[f,g] \mathcal N_{\divi(m)|_{X_K}}(f)\mathcal N_{\divi(l)|_{X_K}}(g)\left<l,m\right>_K$$
therefore by proposition \ref{restr_prop} we can conclude that:
$$\left<l',m'\right>_{U}=\mathcal N_{(f)|_{X_{U}}}(g) \mathcal N_{\divi(m)|_{X_{U}}}(f)\mathcal N_{\divi(l)|_{X_{U}}}(g)\left<l,m\right>_U=a\left<l,m\right>_{U}\,.$$
\endproof

We are ready to show that $\left<\mathscr L,\mathscr M\right>$ is an invertible sheaf on $S$. By  proposition \ref{sstep1} $\left<\mathscr L,\mathscr M\right>$ is non-zero; now assume $\mathscr L=\mathscr O_X(D)$, $\mathscr M=\mathscr O_X(E)$ and fix a point $s_0\in S$. By  the moving lemma we can find a divisor $D'$ such that $D'\sim D$ and $D'$ doesn't have components in $X_{s_0}$. Suppose that $x_1,\ldots x_m$ are the intersection points of $D'$ and $X_{s_0}$, by applying again the moving lemma we can find a divisor $E'$ such that: $E'\sim E$, $E'$ and $D'+X_{s_0}$ have no common components, and $E$ doesn't pass by  $x_1,\ldots, x_m$. Consider the finite  subset of $S$ 
$$C:=\{s\in S'\colon D'\cap E'\cap X_s\neq \emptyset\}\,$$
and note that its complement $U:=S\setminus C$ has the following properties: $s_0\in U$ and $\left<D'|_U, E'|_U \right>=0$. At this point any two meromorphic sections of $\mathscr L$ and $\mathscr M$ corresponding respectively to the divisors $D'$ and $E'$ will satisfy the hypothesis of proposition \ref{sstep2} on $U$. This implies that $\left<\mathscr L,\mathscr M\right>$ is an invertible sheaf.

One can show that the pairing constructed above satisfies the axioms (1)-(4) of definition \ref{axiom_deligne} for $n=1$, moreover we have the following additional properties:

\begin{theorem}\label{del_thm}
The Deligne pairing $(\mathscr L,\mathscr M)\to \left<\mathscr L,\mathscr M\right>$ satisfies the properties listed below. We assume that $\mathscr L$ and $\mathscr M$ are two invertible sheaves on $X$.
\begin{itemize}
\item[$(1)$] The induced map  $\Pic(X)\times\Pic(X)\to\Pic(S)$ is bilinear and symmetric. 
\item[$(2)$] Let $l$ and $m$ be two non-zero meromorphic sections  of $\mathscr L$ and $\mathscr M$, respectively, such that $\divi(l)$ and $\divi(m)$ have no common components. Then,  there exists a non-zero meromorphic section $\left<l,m\right>$  with the following properties: 
\begin{itemize}
\item[$(i)$] If $f,g\in K(X)^\times$ such that $(\divi(fl),\divi(m)),(\divi(l),\divi(gm))\in \Upsilon$, then:
$$\left<fl,m\right>=\mathcal N_{\divi(m)}(f)\left<l,m\right>$$
$$\left<l,gm\right>=\mathcal N_{\divi(l)}(g)\left<l,m\right>$$
\item[$(ii)$] There is an isomorphism of invertible sheaves 
$$\left<\mathscr L,\mathscr M\right>\cong\mathscr O_S(\left<\divi(l),\divi(m)\right>)$$
Moreover, under the above isomorphism $\left<l,m\right>$ corresponds to $1_{\left<\divi(l),\divi(m)\right>}$. In particular:
$$\divi(\left<l,m\right>)=\left<\divi(l),\divi(m)\right>\,.$$
\end{itemize}
\end{itemize}
\end{theorem}
\proof
See \cite[Theorem 4.7]{Mor}.
\endproof
\begin{remark}\label{del_vectspace}
Note that when $S=\spec k$ for any field $k$ (in other words $X$ is an algebraic curve), then $\left<\mathscr L,\mathscr M\right>$ is just a one dimensional $k$-vector space.
\end{remark}
 
\begin{remark}
If $S$ is a non-singular projective curve over a field $k$, and $\varphi: X\to S$ is a morphism over $\spec k$ (i.e. $X$ is a fibred surface over $S$), then it is evident that $(D,E)\mapsto \deg\left<\mathscr O_X(D),\mathscr O_X(D)\right>$ satisfies all the axioms of definition \ref{inters_var}, so it is the intersection pairing on $X$. The same argument holds also when $X$ has generic dimension $n$.
\end{remark}

\end{appendices}

\bibliographystyle{hplain.bst}
\bibliography{expl_del.bib}

\begin{thebibliography}{10}

\bibitem{Beau}
A.~Beauville.
\newblock {\em Complex algebraic surfaces}.
\newblock Cambridge University Press, Cambridge New York, 1996.

\bibitem{SGA}
P.~Berthelot, A.~Grothendieck, L.~Illusie, et~al.
\newblock {\em Th\'{e}orie des intersections et th\'{e}or\`{e}me de
  Riemann-Roch}.
\newblock Springer-Verlag, Berlin, New York, 1971.

\bibitem{BSW}
I.~Biswas, G.~Schumacher, and L.~Weng.
\newblock Deligne pairing and determinant bundle.
\newblock {\em Electron. Res. Announc. Math. Sci.}, 18:91--96, 2011.

\bibitem{bouck}
S.~Boucksom and D.~Eriksson.
\newblock Spaces of norms, determinant of cohomology and {F}ekete points in
  non-{A}rchimedean geometry.
\newblock {\em Adv. Math.}, 378, 2021.

\bibitem{Ca}
P.~Cartier.
\newblock Sur un th\'{e}or\`eme de {S}napper.
\newblock {\em Bull. Soc. Math. France}, 88:333--343, 1960.

\bibitem{Del}
P.~Deligne.
\newblock Le d\'{e}terminant de la cohomologie.
\newblock In {\em Current trends in arithmetical algebraic geometry ({A}rcata,
  {C}alif., 1985)}, volume~67 of {\em Contemp. Math.}, pages 93--177. Amer.
  Math. Soc., Providence, RI, 1987.

\bibitem{Ducr}
F.~Ducrot.
\newblock Cube structures and intersection bundles.
\newblock {\em J. Pure Appl. Algebra}, 195(1):33--73, 2005.

\bibitem{Elk}
R.~Elkik.
\newblock Fibr\'{e}s d'intersections et int\'{e}grales de classes de {C}hern.
\newblock {\em Ann. Sci. \'{E}cole Norm. Sup. (4)}, 22(2):195--226, 1989.

\bibitem{Gar}
E.~Mu{\~n}oz Garcia.
\newblock Fibr\'{e}s d'intersection.
\newblock {\em Compositio Mathematica}, 124(3):219–252, 2000.

\bibitem{EGAIV}
A.~Grothendieck.
\newblock {\'E}l\'ements de g\'eom\'etrie alg\'ebrique : {IV}. \'etude locale
  des sch\'emas et des morphismes de sch\'emas, quatri\`eme partie.
\newblock {\em Publications Math\'ematiques de l'IH\'ES}, 32:5--361, 1967.

\bibitem{Hoche}
A.~Hochenegger.
\newblock Appendix: Introduction to derived categories of coherent sheaves.
\newblock {\em Birational Geometry of Hypersurfaces}, pages 267--295, 2019.

\bibitem{Kl}
S.~L. Kleiman.
\newblock {T}he {P}icard {S}cheme.
\newblock 2014, 1402.0409.

\bibitem{KM}
F.~F. Knudsen and D.~Mumford.
\newblock The projectivity of the moduli space of stable curves. {I}.
  {P}reliminaries on ``det'' and ``{D}iv''.
\newblock {\em Math. Scand.}, 39(1):19--55, 1976.

\bibitem{Kol}
J.~Koll\'{a}r.
\newblock {\em Rational curves on algebraic varieties}.
\newblock Springer, Berlin New York, 1996.

\bibitem{Liu}
Q.~Liu.
\newblock {\em Algebraic Geometry and Arithmetic Curves (Oxford Graduate Texts
  in Mathematics)}.
\newblock Oxford University Press, 8 2006.

\bibitem{MB}
L.~Moret-Bailly.
\newblock M\'{e}triques permises.
\newblock Number 127, pages 29--87. 1985.
\newblock Seminar on arithmetic bundles: the Mordell conjecture (Paris,
  1983/84).

\bibitem{Mor}
A.~Moriwaki.
\newblock {\em Arakelov Geometry (Translations of Mathematical Monographs)}.
\newblock American Mathematical Society, 11 2014.

\bibitem{Phon}
D.~H. Phong, J.~Ross, and J.~Sturm.
\newblock Deligne pairings and the {K}nudsen-{M}umford expansion.
\newblock {\em J. Differential Geom.}, 78(3):475--496, 2008.

\bibitem{Sn1}
E.~Snapper.
\newblock Multiples of divisors.
\newblock {\em J. Math. Mech.}, 8:967--992, 1959.

\bibitem{Sn2}
E.~Snapper.
\newblock Polynomials associated with divisors.
\newblock {\em J. Math. Mech.}, 9:123--139, 1960.

\bibitem{Xia}
M.~Xia.
\newblock {D}eligne-{R}iemann-{R}och {T}heorems {I}. {U}niqueness of {D}eligne
  {P}airings and {D}egree $1$ {P}art of {D}eligne-{R}iemann-{R}och
  isomorphisms, 2017, 1710.09731.

\bibitem{Zh}
S.~Zhang.
\newblock Heights and reductions of semi-stable varieties.
\newblock {\em Compositio Math.}, 104(1):77--105, 1996.

\end{thebibliography}

\end{document}